\documentclass[a4paper,11pt]{amsart}

\usepackage[a4paper, left=2.3cm, right=2.3cm, top=2cm, bottom=2cm]{geometry}
\usepackage{meintex}
\usepackage{caption}
\usepackage{subcaption}
\usepackage[algo2e,norelsize]{algorithm2e}
\usepackage{algpseudocode}
\usepackage[section]{placeins}

\title[On the Numerics of Mean-field Optimal Control]{Controlling a self-organizing system of individuals  guided by a few external agents --- Particle  description and mean-field limit }

\author[M.~Burger]{Martin Burger}
\address[M.~Burger]{Westf\"alische Wilhelms-Universit\"at M\"unster}
\email{martin.burger@wwu.de}

\author[R.~Pinnau]{Ren\'e Pinnau}
\address[R.~Pinnau]{TU Kaiserslautern}
\email{pinnau@mathematik.uni-kl.de}

\author[A.~Roth]{Andreas Roth}
\address[A.~Roth]{TU Kaiserslautern}
\email{roth@mathematik.uni-kl.de}

\author[C.~Totzeck]{Claudia Totzeck}
\address[C.~Totzeck]{TU Kaiserslautern}
\email[Corresponding author]{totzeck@mathematik.uni-kl.de}

\author[O.~Tse]{Oliver Tse}
\address[O.~Tse]{TU Kaiserslautern}
\email{tse@mathematik.uni-kl.de}

\begin{document}
\maketitle
\begin{abstract}
Optimal control of large particle systems with collective dynamics by few agents is a subject of high practical importance (e.g. in evacuation dynamics), but still limited mathematical basis. In particular the transition from discrete optimal control to a continuum setting as the number of particles tends to infinity is by far not fully understood. In this paper we contribute to this issue by studying a canonical model of controlling an interacting particle system into a certain spatial region by repulsive forces from few external agents, which might be interpreted as shepherd dogs leading sheep to their home. 

We discuss the appropriate modelling of such a problem and the associated optimality systems, providing some connections between the Lagrange multipliers in the discrete and continuum setting. As control strategies we investigate an {\em Instantaneous Control} and a global {\em Optimal Control} approach.
 The solutions of a family of control problems for the particle system with external agents are numerically compared to the mean-field controls as the number of particles tends to infinity. In both cases, this leads to a high dimensional phase space requiring tailored optimization strategies.  All control problems arising are solved using adjoint information to compute the descent directions. The numerical results indicate the convergence of controls for both optimization strategies.
\end{abstract}
\section{Introduction}
In the last decades, the behavior of large particle systems and their mean-field limits were extensively investigated with theoretical and computational approaches \cite{Golse,AlbiPareschi,CarrilloMills,carrillo2009double}. 
Large groups of individuals like flocks of birds and schools of fish, and their attractive and repulsive interaction were considered, which lead to models of different types of collective behavior such as flocking or milling, and a thorough study of their stability \cite{reynolds1987flocks,czirok2000collective,barbaro2009discrete,maini2001mathematical} (also \cite{CarrilloMills} for detailed a overview).
The models were refined using vision cones, self-propulsion and orientation alignment of neighbors\cite{d2006self,CarrilloMills,yates2011refining}. As the behavior of large particle groups of same type is well understood, the interest in self-organized systems interacting with few external agents arose. This concept and its kinetic limit was first numerically investigated by Albi and Pareschi in \cite{AlbiPareschi}. Using interaction potentials introduced by Cucker-Smale \cite{cucker2007emergent} or D'Orsogna et al \cite{Dorsogna}, they showed numerically, that the collective behavior of large groups coincides with the behavior of the kinetic model as the number of individuals tends to infinity. While the models for the interaction of particles became more and more realistic, the complexity was significantly increased. At first, interaction potentials were chosen smooth in order to have well-defined derivatives \cite{cucker2007emergent,cucker2007mathematics,Dorsogna}. Later it was shown that special classes of singular interaction potentials allow to pass to the mean-field limit as well \cite{hauray2007n,bolley2011stochastic}. 

In the present work, we investigate different control strategies for the interaction of a huge crowd of individuals with few external agents, with a focus on modelling and computation. This behavior is then compared to the corresponding mean-field  control problem. As control parameters we use the velocities of the external agents. The cost functionals are designed in such a way that the external agents aim to lead the crowd to a predefined destination using the least amount of energy possible. As an example one may consider a crowd of sheep guided by dogs. 
Analytically,  a similar sparse optimization problem was investigated by Fornasier and Solombrino in \cite{FornasierSolombrino}. They showed the existence of optimal controls on the microscopic level and used the concept of $\Gamma$-convergence to perform the mean-field limit in a sparse optimal control setting.

Numerically, already the computation of the state solutions are challenging due to the high dimensional phase space. On the one hand, we have the pairwise interactions in the microscopic system and on the other hand the mean-field equations are of Vlasov type, which yields in two spatial dimensions already a four dimensional problem. This becomes even worse for the optimization of such problems, since the state systems need to be solved several times and a huge amount of data needs to be stored. Hence, there is no chance to solve these problems using black-box optimization approaches. Instead, tailored optimization algorithms are required, which are typically based on derivative information \cite{Pinnau}. 

In the literature, different control strategies are available to investigate such problems. For example the Instantaneous Control (IC) was applied to traffic flow problems \cite{HertyKlar} and the Navier-Stokes Equations \cite{HinzeNavier}. The Optimal Control (OC) approach was used in modeling of glass cooling \cite{Pinnau} or semiconductors (cf. \cite{DebyeGamma} and the overview in \cite{Semiconductor}). Due to the large amount of memory needed to store the forward information, the OC approach is not frequently used for controlling multi-dimensional problems as considered here.

The major aim of this manuscript is to construct appropriate numerical approaches and computationally verify that the controls of the microscopic problem converge to the optimal control of the kinetic problem as the number of individuals increases for each of the strategies. To this end we implemented numerical descent algorithms, based on first-order derivative information. These are derived by applying the adjoint calculus to the interacting particle system as well as to the mean-field equation.  A simple application scenario we have in mind is to guide a herd of sheep into a given region (stable) by few dogs, the interaction between the particles and the external agents thus being of repulsive type.

The manuscript is organized as follows. The state systems which constrain the optimization problems are defined for the microscopic and the kinetic case in Section 2. The cost functionals and the corresponding constrained optimization problems are stated in Section 3.  In Section 4 the associated first order optimality conditions are derived. Section 5 addresses the numerical schemes and the optimization algorithms. The numerical results for the IC and the OC approach are discussed in Section 6. Concluding Remarks are given in Section 7.
\section{Microscopic and Mean-Field Optimal Control Problems}\label{sec:state}
We start by describing the agent-based model and its corresponding mean-field limit, which are the state systems for the control problems considered later on.

\subsection{Microscopic Model}
Let $D\ge 1$ denote the dimension of the spatial and velocity space. The considered particle system consists of $N \in \mathbb{N}$ particles of the same type and $M$ external agents with $N\gg M$. Let $[0,T]$ be the time interval of observation, then the particles and agents are represented by state vectors
\begin{equation*}
 x_i,v_i,d_m,u_m \colon [0,T] \rightarrow \mathbb{R}^D, \qquad \text{for } i = 1,\dots N \text{ and }m=1,\dots,M.
\end{equation*}
The vectors 
\begin{equation*}
 \x(t) = (x_i(t))_{i=1,\dots,N}\quad \text{ and } \quad \ve(t) = (v_i(t))_{i=1,\dots,N}, 
\end{equation*}
 denote the positions and velocities of the particles and 
 \begin{equation*}
 \dog(t) = (d_m(t))_{m=1,\dots,M}\quad \text{ and } \quad \con(t) = (u_m(t))_{m=1,\dots,M},
 \end{equation*}
 the positions and velocities of the external agents, respectively (see also \cite{AlbiPareschi}). Since the time dependence is clear, we often write $\x$ or 
$x_i$ instead of $\x(t)$ respectively $x_i(t)$ to shorten the notation. Note that we write the vectors $\x$ and $\ve$ in bold to make clear when denoting the vectors containing positions and velocities of individuals and when we refer to the position and velocity space with variables $x$ and $v$ in the mean-field setting later on.

The interactions of the individuals and external agents are modeled using potentials $\Phi_j$, $j=1,2$ which satisfy the following assumption
\begin{itemize} 
 \item [\textbf{(A)}] $\Phi_j\colon\mathbb{R}^D\to\mathbb{R}$ are radially symmetric and continuously differentiable, with $\nabla \Phi_j$ locally Lipschitz and globally bounded.
\end{itemize}
These rather strict assumptions allow for existence and uniqueness of solutions to the adjoint and state systems. The latter are essential for defining the reduced cost functional needed in the algorithms for numerical investigation. In addition to the interaction terms, our model includes (linear) friction represented by the friction parameter $\alpha>0$. To simplify the presentation, we denote for any $x,y\in\mathbb{R}^D$ the interaction forces 
\[
 K_j(x,y) = (\nabla\Phi_j)(x-y),\qquad j=1,2.
\]
Since $\Phi_j$ is symmetric, we have that $K_j(x,y)=-K_j(y,x)$ for $x,y\in\mathbb{R}^D$.

Altogether, this leads to the particle system
\begin{subequations}\label{longODEsys}
 \begin{align}
 \frac{\dd}{\dd t}x_i &= v_i,\qquad \frac{\dd}{\dd t}v_i = - \frac{1}{N} \sum_{k\ne i} K_1(x_i,x_k) - \frac{1}{M} \sum_{m=1}^M K_2(x_i,d_m) - \alpha v_i, \label{ODExv}\\
 \frac{\dd}{\dd t}d_m &= u_m. \label{ODEd}
\end{align}
\end{subequations}
The individuals interact pairwise $(j=1)$ with each other and with the external agents $(j=2)$. The latter do not interact among themselves nor are influenced by the others. Note however that such interactions can easily be considered by a simple transformation of the optimal control variables analogous to the approach in \cite{burger2009globally}. For system \eqref{longODEsys} the velocities of the external agents $u$ are assumed to be given. Later, they serve as controls for the optimization problem. We assume that the external agents have a  maximal speed $u_{\text{max}}$ which induces a constraint in the space of admissible controls (see \eqref{eq:control} below).

For notational convenience we further define the state vector $y := (\x,\ve,\dog) \in \mathbb{R}^{D(2N+M)}$ and
\begin{equation*}
  {\bf S}(y) = ({\bf S}_i(y))_{i=1,\dots,N},\qquad  {\bf S}_i(y) = -\frac{1}{N} \sum_{j\ne i} K_1(x_i,x_j) - \frac{1}{M} \sum_{m=1}^M K_2(x_i,d_m) - \alpha v_i. 
\end{equation*}
Using this notation, the ODE system \eqref{longODEsys} 
may be written compactly as
\begin{subequations}\label{ODEsys}
\begin{equation}\label{stateODE}
\frac{\dd}{\dd t}{y}= (\frac{\dd}{\dd t}{\x}, \frac{\dd}{\dd t}{\ve}, \frac{\dd}{\dd t}{\dog} )  = (\ve, {\bf S}(y), \con ) =: F(y,\con), 
\end{equation}
with the mapping $F: \mathbb{R}^{D(2N+M)} \times \mathbb{R}^{DM} \to \mathbb{R}^{D(2N+M)}$. The microscopic state system is supplemented with the initial conditions
\begin{equation}\label{IC_ODE}
 \x(0) = \x_0 \in \mathbb{R}^{DN}, \qquad \ve(0) = \ve_0 \in \mathbb{R}^{DN},\qquad \dog(0)=\dog_0 \in \mathbb{R}^{DM}.
\end{equation}
\end{subequations}
These initial conditions are denoted in short by $y(0) = y_0$.

\subsubsection{Remarks on Well-Posedness}
The well-posedness of \eqref{ODEsys} is rather standard. However, since it is common in optimal control problems to consider a Hilbert space for the control parameters, we use
\begin{align}\label{eq:space_control}
 U:=L^2((0,T),\mathbb{R}^{MD})
\end{align}
as space of controls. Nevertheless, this is sufficient to prove the well-posedness of (\ref{ODExv}). Indeed, $\dog$ may be explicitly expressed as
\[
 \dog(t) = \dog_0 + \int_0^t \con(s)\dd s,
\]
which shows that $\dog$ is absolutely continuous and hence also continuous. Consequently, we obtain from the boundedness of $\con$ the existence and uniqueness of a global solution due to the theorem by Picard and Lindel\"of. From the expression above we easily deduce that $\dog\in H^1((0,T),\mathbb{R}^{MD})$.

Altogether, we are able to define the corresponding \textit{control-to-state} operator $\mathcal{G}_N \colon U \rightarrow Y$ which maps any control parameter $\con\in U$ to the unique solution $y=\mathcal{G}_N(\con)$ of \eqref{ODEsys} in the state space $Y$. In this case, we choose the state space
\begin{align}\label{eq:space_state_N}
 Y := H^1((0,T),\mathbb{R}^{ND})\times H^1((0,T),\mathbb{R}^{ND}) \times H^1((0,T),\mathbb{R}^{MD}).
\end{align}
Note, that the solution $y\in \mathcal{C}^1((0,T),\mathbb{R}^{ND})\times \mathcal{C}^1((0,T),\mathbb{R}^{ND}) \times H^1((0,T),\mathbb{R}^{MD})\subset Y$ is indeed in a subspace of $Y$.

\subsection{Mean-Field Model}
In order to define the limiting problem for an increasing number of individuals $N$ explicitly, we consider the empirical measure
\[
 \mu_t^N(x,v) = \frac{1}{N} \sum_{i=1}^N \delta_0(x_i(t)-x) \otimes \delta_0(v_i(t)-v).
\]
By definition, $\mu_t^N\in\mathcal{P}(\mathbb{R}^{2D})$ is a Borel probability measure that assigns the probability $\mu_t^N(A)$ of finding particles with states within a Borel measurable set $A\subset\mathbb{R}^{2D}$ in the phase space $\mathbb{R}^{2D}$ at time $t\ge 0$. If a Borel probability measure $\mu_t\in \mathcal{P}^{ac}(\mathbb{R}^{2D})$ is absolutely continuous w.r.t.~the Lebesgue measure, we denote its density by $f_t\in L^1(\mathbb{R}^{2D})$. For later use we introduce the macroscopic density $\rho_t$ of a Borel probability measure $\mu_t\in\mathcal{P}(\mathbb{R}^{2D})$ as its first marginal, i.e.,
\begin{equation*}
\rho_t(A) := \mu_t(A\times \mathbb{R}^D) = \iint_{A\times \mathbb{R}^D} f_t(x,v)\dd x\dd v,
\end{equation*}
for any Borel measurable set $A\subset\mathbb{R}^D$. The last equality holds whenever $\mu_t\in\mathcal{P}^{ac}(\mathbb{R}^{2D})$.

The link between the particle system and the mean-field equation is derived formally using ideas from \cite{Dobrushin,BraunHepp,Neunzert}. Let therefore $h \in \mathcal{C}^\infty_c(\mathbb{R}^{2D})$ be an arbitrary real infinitely differentiable function on $\mathbb{R}^{2D}$ with compact support and $z_i=(x_i,v_i)\in \mathbb{R}^{2D}$. Then
\begin{equation*}
 \frac{\dd}{\dd t} \left< \mu_t^N , h \right> := \frac{\dd}{\dd  t} \frac{1}{N}\sum_{i=1}^N h(z_i) = \frac{1}{N}\sum_{i=1}^N \nabla_{x} h(z_i)\cdot\frac{\dd}{\dd t}x_i + \nabla_{v} h(z_i)\cdot\frac{\dd}{\dd t}v_i,
\end{equation*}
which allows for the formal calculation
\begin{align}\label{eq:FNweakSol}
 \langle \partial_t \mu_t^N , h \rangle 
 &= \frac{1}{N} \sum_{i=1}^N  \nabla_{x} h(z_i) \cdot v_i - \nabla_{v} h(z_i) \cdot \left( \big(K_1\ast  \rho_t^N\big)(x_i) + \frac{1}{M} \sum_{m=1}^M K_2(x_i,d_m) +\alpha v_i\right) \nonumber\\
 &= \left< \mu_t^N , \nabla_x h \cdot v - \nabla_v h \cdot \left( \big(K_1\ast  \rho_t^N\big)(x) + \frac{1}{M} \sum_{m=1}^M K_2(x,d_m) +\alpha v\right) \right>,
 \end{align}
 where 
 \[
  \big(K_1\ast  \rho_t^N\big)(x) = \int_{\mathbb{R}^{2D}} K_1(x,\bar{x}) \dd \mu_t^N(\bar{z}).
 \]
Passing to the limit $N\to\infty$ and integrating by parts, we arrive at the equation
 \[
 0 = \langle \partial_t \mu_t + v \cdot \nabla_x \mu_t + \nabla_v \cdot (S(\mu_t)  \mu_t) , h\rangle,
 \] 
 where we define
 \[
  S(\mu_t)(x,v,d) = -\big(K_1\ast  \rho_t\big)(x) - \frac{1}{M} \sum_{m=1}^M K_2(x,d_m) - \alpha v.
 \]
Since $h \in \mathcal{C}^\infty_c(\mathbb{R}^{2D})$ is arbitrary, we may use the variational lemma to find that
\begin{equation}\label{longPDE}
 \partial_t \mu_t  + v\cdot \nabla_x \mu_t + \nabla_v  \cdot (S(\mu_t) \mu_t )=0,
\end{equation}
which is the mean-field single particle distribution, 
supplemented with the initial condition $\mu(0)=\mu_0$.

\begin{rem}
Observe that $\mu_t^N$ and $\mu_t$ satisfy exactly the same equation in the distributional sense. The above limit exists for example in the Wasserstein metric (for a detailed discussion see \cite{WassersteinConvergence}). Note that even though $\mu_t$ has unit mass by construction, we do not need to explicitly assume that the condition is satisfied, since this is a natural constraint and the control parameter does not change this property. For this reason the use of the standard theory of OC based on the $L_2$-calculus is justified. Otherwise we would enter the regime of probability measures leading to Wasserstein derivatives. Note that while we pass to the mean-field limit $N \rightarrow \infty$, the number of external agents $M$ remains finite. 
\end{rem}

%

\subsubsection{Remarks on Well-Posedness} The existence and uniqueness of solutions for the Vlasov equation \eqref{longPDE} may be found, for example, in the articles \cite{BraunHepp,WassersteinConvergence,Dobrushin,Golse}, where the notion of solution is established in the Wasserstein space of Borel probability measures. 

\begin{defi}
 Let $\mathcal{P}_1(\mathbb{R}^{2D})$ denote the space of Borel probability measures on $\mathbb{R}^{2D}$ with finite first moment. We say that $\mu\in \mathcal{C}([0,T],\mathcal{P}_1(\mathbb{R}^{2D}))$ is a {\em weak measure solution} of \eqref{longPDE} with initial condition $\mu_0\in \mathcal{P}_1(\mathbb{R}^{2D})$ if for any test function $h\in \mathcal{C}_c^\infty((-\infty, T]\times \mathbb{R}^{2D})$ we have
 \[
  \int_0^T \int_{\mathbb{R}^{2D}} \big( \partial_t h_t + v\cdot\nabla_x h_t + S(\mu_t)\cdot\nabla_v h_t \big)\dd \mu_t\dd t + \int_{\mathbb{R}^{2D}} h_0\dd \mu_0 = 0.
 \]
\end{defi}

On the other hand, the equation \eqref{longPDE} may be equivalently expressed as a nonlinear flow
\begin{align}\label{eq:nonliner_process}
 \frac{\dd}{\dd t} Z = (\frac{\dd}{\dd t} x,\frac{\dd}{\dd t} v) = \big( v, S(\mu_t)(Z,d)\big),\qquad \mu_t=\text{law}(Z(t)),
\end{align}
 with the initial condition $Z(0)=Z_0\in\mathbb{R}^{2D}$, where $Z_0$ is a random variable distributed according to $\mu_0=\text{law}(Z_0)$. Assuming the solvability of the nonlinear process \eqref{eq:nonliner_process}, a weak measure solution of $\eqref{longPDE}$ may be represented as the push-forward of the measure along the flow $Z(t,Z_0)$, i.e., $\mu_t=Z(t,\cdot)\# \mu_0$. If we further assume that the initial measure $\mu_0\in\mathcal{P}_1(\mathbb{R}^{2D})$ has density $f_0\in L^1(\mathbb{R}^{2D})$, and that the nonlinear flow given by \eqref{eq:nonliner_process} satisfies $Z\in \mathcal{C}([0,T],\text{Diff}(\mathbb{R}^{2D}))$, i.e., $Z(t,\cdot)$ is a diffeomorphism for all $t\ge 0$, then $\mu_t\in\mathcal{P}_1^{ac}(\mathbb{R}^{2D})$ has density $f_t\in L^1(\mathbb{R}^{2D})$ and $f\in \mathcal{C}([0,T],L^1(\mathbb{R}^{2D}))$.

In order to employ the standard $L_2$-calculus we will require more regularity of the states and  assume additionally that
\begin{itemize}
 \item [\textbf{(B)}] $\mu_0\in\mathcal{P}_1(\mathbb{R}^{2D})$ has density $f_0\in L^2(\mathbb{R}^{2D})$ with compact support and the flow given by \eqref{eq:nonliner_process} satisfies $Z\in \mathcal{C}^1([0,T],\text{Diff}(\mathbb{R}^{2D}))$ such that $f\in H^1((0,T),L^2(\mathbb{R}^{2D}))$.
\end{itemize}

\begin{rem}
 Note that the assumption on the nonlinear flow $Z\in \mathcal{C}^1([0,T],\text{Diff}(\mathbb{R}^{2D}))$ satisfying \eqref{eq:nonliner_process} in {\bf(B)} is not restrictive due to assumption {\bf(A)}. Indeed, since $\mu\in \mathcal{C}([0,T],\mathcal{P}_1(\mathbb{R}^{2D}))$, we have that $S(\mu_t)$ is continuous in $t$ and locally Lipschitz in $z$. Therefore, standard ODE theory provides the required regularity for the nonlinear flow.
\end{rem}

Defining the state variable $p:= (f,\dog)$, the coupled system of \eqref{longPDE} and the ODE \eqref{ODEd} that models the movement of the external agents $\dog$ with control ${\bf w}$ may be written as
\begin{subequations}\label{PDEsys}
\begin{equation} \label{statePDE}
 \frac{\dd}{\dd t} p := (\partial_t f_t , \frac{\dd}{\dd t}\dog ) = -\big(\nabla_v \cdot (S(f_t)  f_t) + v \cdot \nabla_x f_t, {\bf w}\big) =: G(p,{\bf w}).
\end{equation}
The initial conditions for this system are
\begin{equation}\label{IC_PDE}
f_t|_{t=0} = f_0\in L^2(\mathbb{R}^{2D}), \qquad \dog(0) = \dog_0\in\mathbb{R}^{MD}. 
\end{equation}
\end{subequations}
In the following we shall refer to \eqref{IC_PDE} as $p(0)=p_0$. 

Since the problem is defined on the whole space $\mathbb{R}^{2D}$, there are no boundary conditions required. In fact, from {\bf(B)}, we expect $f_t$ to have compact support that remains within a bounded convex domain $\Omega\subset \mathbb{R}^{2D}$ with smooth boundary for all times $t\in[0,T]$ (cf.~\cite[Thm 4.4]{FornasierSolombrino}). For this reason, we may define the corresponding \textit{control-to-state} operator $\mathcal{G}_\infty\colon U\to\mathcal{Y}$, $\con\mapsto p=\mathcal{G}_\infty(\con)$, where $p$ satisfies \eqref{PDEsys} with the state space $\mathcal{Y}$ given by
\begin{align}\label{eq:state_vlasov}
 \mathcal{Y}:= H^1((0,T),L^2(\Omega))\times H^1((0,T),\mathbb{R}^{MD}).
\end{align}
Here, $\Omega\subset\mathbb{R}^{2D}$ is an a priori given bounded convex domain with smooth boundary which contains the support of $f_t$ for all times $t\in[0,T]$.  


%
%

\section{The Optimal Control Problems}

The aim of the optimal control problem is to guide the crowd to a specified destination $\Edes \in \mathbb{R}^D$ by controlling the velocity of the external agents. In order to model this issue mathematically as an optimal control problem, we will need a cost functional that reflects the aim of steering the system of particles to a desired destination on the one hand, and on the other hand it has to converge in the limit $N \rightarrow \infty$ to a suitable cost functional on the mean-field level. In the following, we state one possible cost functional based on the center of mass and the variance of the crowd, which meets these conditions.

\subsection{Optimal Control Problem in the Microscopic Setting}
As above $\x$ denotes the particle positions. Hence the center of mass $\mathbb{E}_N$ and variance $\mathbb{V}_N$ for the particles are defined as 
\begin{equation*}
 \mathbb{E}_N(\x(t)) = \frac{1}{N}  \sum_{i=1}^N x_i(t),\qquad \mathbb{V}_N(\x(t)) = \frac{1}{N} \sum_{i=1}^N |x_i(t)-\mathbb{E}_N(\x(t))|^2,
\end{equation*}
and we define the cost functional
\begin{equation}\label{eq:cost functional_N}
 J_N(y(\con),\con) = \frac{1}{T} \int_0^T \frac{\sigma_1}{4} |\mathbb{V}_N(\x(t)) - \bar{V}_N|^2+ \frac{\sigma_2}{2}  |\mathbb{E}_N(\x(t))-\Edes|^2  + \frac{\sigma_3}{2M} \norm{\con(t)}_{\mathbb{R}^{MD}}^2 \dd t,
\end{equation}
where $\bar{V}_N$ is the desired variance, e.g. zero. The first part of the functional penalizes the difference of the spread of the particles with respect to the desired variance $\bar{V}_N$. The second term measures the distance of the center of mass to the destination $\Edes\in\mathbb{R}^D$. Technically, the third part ensures the local convexity of the cost functional. Physically, it is an energy term causing the velocities of the external agents to remain bounded. Thus, the energy supplied to the system by the movement of the agents is minimized. 

Altogether, the above cost functional models the task of guiding the crowd of individuals such that the crowd is clustered around $\Edes\in\mathbb{R}^D$ while utilizing the least amount of kinetic energy possible. The positive parameters $\sigma_i$ allow to adjust the influence of the different parts of the cost functional. Due to the maximal speed constraint on the velocities of the external agents we define the space of admissible controls $\mathcal{U}_{ad}$ as 
\begin{equation}\label{eq:control}
 \mathcal{U}_{ad} := \left\{ \con \in L^2((0,T),\mathbb{R}^{MD})\; \colon\; |u_m(t)| \le u_{\text{max}}, \quad m=1,\dots,M \right\}.
\end{equation}
This yields a family of constrained optimal control problems on the microscopic level given by
\begin{problem}\label{Opt_ODE}
 Find $(y_*,\con_*) \in Y\times \mathcal{U}_{ad}$ such that
 \begin{equation*} 
  (y_*,\con_*) = \argmin\nolimits_{(y,\con) \in Y\times \mathcal{U}_{ad}}  J_N(y,\con)\quad \text{ subject to the IVP } \eqref{ODEsys}, \\
 \end{equation*}
 where $J_N(y,\con)$ is given by \eqref{eq:cost functional_N}.
\end{problem}

\begin{rem}
 For the IC approach we cut the time interval into slices and solve a stationary control problem on each time slice (see Section \ref{sec:Alg}). The velocities of the external agents are assumed to be constant on each time slice. The function $\con \in \mathcal{U}_{ad}$ is the piecewise constant collection of the solutions on the slices projected onto the feasible set. Due to the constraint on the maximal velocity  \eqref{eq:control}, we use a projected gradient method to solve the optimization problems. 
\end{rem}

\subsection{Optimal Control Problem in the Mean-Field Setting}
As discussed already at the beginning of this section, the control problem on the microscopic level has to match the one of the mean-field level as $N \rightarrow \infty$. Thus, the following cost functional is a direct consequence of \eqref{eq:cost functional_N}. The center of mass $\mathbb{E}_\infty$ and variance $\mathbb{V}_\infty$ only depend on the macroscopic density and are defined as 
\begin{equation*}
 \mathbb{E}_{\infty}(f_{t}) = \iint_{\mathbb{R}^{D}\times \mathbb{R}^D} x  f_t\dd x \dd v,\qquad \mathbb{V}_{\infty}(f_{t}) = \iint_{\mathbb{R}^{D} \times \mathbb{R}^D} |x -\mathbb{E}_{\infty}(f_{t})|^2 f_t\dd x \dd v,
\end{equation*}
leading to the cost functional
\begin{equation}\label{eq:cost functional_MF}
 J_\infty(p({\bf w}),{\bf w}) = \frac{1}{T} \int_0^T \frac{\sigma_1}{4} |\mathbb{V}_{\infty}(f_t) - \bar{V}_\infty |^2+ \frac{\sigma_2}{2}  |\mathbb{E}_{\infty}(f_t)-\Edes|^2 + \frac{\sigma_3}{2M} \norm{\textbf{w}_t}_{\mathbb{R}^{MD}}^2 \dd t.
 \end{equation}
Finally, we define the control problem in the mean-field setting as
\begin{problem}\label{Opt_PDE}
 Find $(p_*,{\bf w}_*) \in \mathcal{Y}\times\mathcal{U}_{ad}$ such that
 \begin{equation*} 
  (p_*,{\bf w}_*) = \argmin\nolimits_{(p,{\bf w}) \in \mathcal{Y}\times \mathcal{U}_{ad}}  J_\infty(p,{\bf w})\quad \text{ subject to the mean-field system \eqref{PDEsys}},\\
 \end{equation*}
 where $J_\infty(p,{\bf w})$ is defined in \eqref{eq:cost functional_MF}.
 \end{problem}

\begin{rem}
An existence result for Problem \ref{Opt_ODE} and \ref{Opt_PDE} may be deduced in a straight-forward way from the results in \cite{FornasierSolombrino}. Even though the authors work in a sparse control setting, similar arguments may be applied in the present setting as well.
\end{rem}
\section{First Order Necessary Conditions} \label{sec:adjoint}
As mentioned in the introduction we apply adjoint based descent methods to solve the control problems. In this section the adjoints and optimality conditions for the particle and the mean-field optimal control problem are formally derived with the  help of the extended Lagrangians corresponding to Problem \ref{Opt_ODE} and \ref{Opt_PDE}. Furthermore, we introduce the reduced cost functional and its gradient.

\subsection{Adjoint of the Microscopic System}
We begin by recalling the spaces $U$ and $Y$ defined in Section~\ref{sec:state}. Let the control space $U$ and state space $Y$ be the Hilbert spaces
\begin{equation*} 
 U=L^2((0,T),\mathbb{R}^{MD}),\qquad Y = [H^1((0,T),\mathbb{R}^{ND})]^2 \times H^1((0,T),\mathbb{R}^{MD}),
\end{equation*}
with $\mathcal{U}_{ad}\subset U$ defined in \eqref{eq:control}. We further denote $X:=[L^2((0,T),\mathbb{R}^{ND})]^2\times L^2((0,T),\mathbb{R}^{MD})$ and
\[
 Z:=X\times\big([\mathbb{R}^{ND}]^2\times \mathbb{R}^{MD}\big),
\]
as the space of Lagrange multipliers with $Z^*$ being its dual. 

We define the state operator $e_N\colon Y \times U \rightarrow Z^*$ of the microscopic problem as
\begin{equation*}
e_N(y,\con) = \begin{pmatrix} \frac{\dd}{\dd t}y - F(y,\con) \\ y(0)-y_0 \end{pmatrix}
\end{equation*}
and the dual pairing
\begin{equation*}
 \langle e_N(y,\con),(\xi,\eta) \rangle_{Z^*,Z} = \int_0^T (\frac{\dd}{\dd t}y(t) - F(y(t),\con(t))) \cdot \xi(t) \dd t + (y(0) - y_0) \cdot \eta.
\end{equation*}
Let $(\xi,\eta)\in Z$ denote the Lagrange multiplier which is in fact the adjoint state. Then, the extended Lagrangian corresponding to Problem \ref{Opt_ODE} reads
\begin{equation*}
 \mathcal{L}_{N}(y,\con,\xi,\eta) = J_{N}(y,\con) + \langle e_N(y,\con),(\xi,\eta) \rangle_{Z^*,Z}.
\end{equation*}
As usual the first-order optimality condition of Problem \ref{Opt_ODE} is derived by solving
\begin{equation*}
 d \mathcal{L}_N(y,\con,\xi,\eta) =0. 
\end{equation*}
The derivative w.r.t.~the adjoint state results in the state equation while the derivative with respect to the state $y$ yields the adjoint system and the optimality condition is obtained by the derivative w.r.t.~the control $\con$ \cite{Pinnau}. 

For the calculations we denote the three parts of the cost functional by $J^i$, $i=1,2,3$, as
\begin{gather*}
 J_N^1(y) = \frac{\sigma_1}{4T} \int_0^T |\mathbb{V}_N(\x(t)) - \bar{V}_N|^2 \dd t, \qquad J_N^2(y) =  \frac{\sigma_2}{2T} \int_0^T \norm{\mathbb{E}_N(\x(t))-\Edes}_{\mathbb{R}^{D}}^2 \dd t,\\
  J_N^3(\con)=\frac{\sigma_3}{2MT} \int_0^T \norm{\con(t)}_{\mathbb{R}^{MD}}^2 \dd t.
\end{gather*}

For any $h=(h^y=(h^x,h^v,h^d),h^u)\in Y\times U$, the following G\^ateaux derivatives are obtained
\begin{align*}
 d_{\x} J_N^1(y) [h^x] 
 &=  \frac{\sigma_1 }{NT} \int_0^T  (\mathbb{V}(\x(t)) - \bar{V}_N)(\x(t) - \mathbb{E}(\x(t))) \cdot h^x(t) \dd t,  \\
 d_{\x} J_N^2(y) [h^x] 
 &= \frac{\sigma_2}{NT} \int_0^T (\mathbb{E}(\x(t)) - \Edes) \cdot  h^x(t) \dd t, \\
 d_{\con} J_N^3(\con) [h^u] &=  \frac{\sigma_3}{M} \int_0^T \con(t) \cdot h^u(t) \dd t.
\end{align*}
For the second part of the Lagrangian we derive
\begin{align*}
 \langle d_\x e_N(y,\con)[h^x],(\xi,\eta) \rangle &= \int_0^T \frac{\dd}{\dd t}h^x(t)\cdot\xi_1(t) - d_{\x} {\bf S}(y)[h^x(t)]\cdot \xi_2(t) \dd t + h^x(0)\cdot\eta_1, \\
 \langle d_{\ve} e_N(y,\con)[h^v], (\xi,\eta) \rangle &= \int_0^T \left(\frac{\dd}{\dd t}h^v(t) + \alpha h^v(t)\right)\cdot \xi_2(t) - h^v(t) \cdot \xi_1(t)  \dd t + h^v(0)\cdot\eta_2, \\
 \langle d_{\dog} e_N(y,\con)[h^d], (\xi,\eta) \rangle &= \int_0^T \frac{\dd}{\dd t}h^d(t)\cdot\xi_3(t) - d_{\dog}{\bf S}(y)[h^d(t)] \cdot \xi_2(t) \dd t + h^d(0)\cdot\eta_3, \\
 \langle d_{\con} e_N(y,\con)[h^u] , (\xi,\eta) \rangle &= -\int_0^T h^u(t)\cdot \xi_3(t) \dd t,
\end{align*}
where the vanishing derivatives are omitted. Assuming that $\xi\in Y$ one may formally derive the strong formulation of the adjoint system. Indeed, integrating by parts and using the fact that $d_{\x} {\bf S}(y)$ and $d_{\dog}{\bf S}(y)$ are symmetric matrices, we arrive at the following result.

\begin{prop}\label{ODEKKT}
 The first-order optimality condition corresponding to Problem \ref{Opt_ODE} reads
 \begin{equation}\label{eq:var_ode}
  \int_0^T \left(\frac{\sigma_3}{MT} \con_*(t) - \xi_3(t) \right)\cdot(\con(t)-\con_*(t)) \dd t\ge 0  \qquad \text{for all\; $\con\in \mathcal{U}_{ad}$},
 \end{equation}
 where $\xi=(\xi_1,\xi_2,\xi_3)\in Y$ satisfies the adjoint system given by
 \begin{subequations}\label{ad_ODE}
 \begin{equation}\label{xi}
  \frac{\dd}{\dd t}\xi_1 = -d_{{\emph \x}} {\bf S}(y)[\xi_2]-d_{\emph \x}J_N(t), \qquad \frac{\dd}{\dd t}\xi_2 = \xi_1 - \alpha \xi_2, \qquad  \frac{\dd}{\dd t}\xi_3 = -d_{{\emph \dog}} {\bf S}(y)[\xi_2],
 \end{equation}
 with
  \begin{equation}\label{terminal_ODE}
  d_{\emph \x} J_N(t) = \frac{\sigma_1}{NT}\Big((\mathbb{V}({\emph \x}(t)) - \bar{V}_N)({\emph \x}(t) - \mathbb{E}_N({\emph \x}(t)))\Big) + \frac{\sigma_2}{NT} \Big(\mathbb{E}_N({\emph \x}(t))-\Edes\Big),
 \end{equation}
supplemented with the terminal conditions $\xi_1(T) = 0$, $\xi_2 (T) = 0$, $\xi_3(T) = 0$.

 \end{subequations}
\end{prop}

\subsection{Adjoint of the Mean-Field System}
 Here, we assume that $p = (f,\dog)$ lies within the state space $\mathcal{Y}$ of the PDE optimization problem, where
 \[
  \mathcal{Y} = H^1((0,T),L^2(\Omega)) \times H^1((0,T),\mathbb{R}^{MD}).
 \]
 Let $\mathcal{X}:=H^1((0,T),L^2(\Omega))\cap L^2((0,T),H^1(\Omega))\times L^2((0,T),\mathbb{R}^{MD})$ and set $\mathcal{Z} := \mathcal{X}\times \left( L^2(\Omega) \times \mathbb{R}^{MD} \right)$ to be the space of adjoint states with dual $\mathcal{Z}^*$. Note that the control space $U$ is identical to the one of the particle problem. 
 
 Now define the mapping $e_\infty\colon \mathcal{Y} \times U \rightarrow \mathcal{Z}^*$ by
\begin{align*}
 \langle e_\infty(p,{\bf w}),(\varphi,\eta) \rangle_{\mathcal{Z}^*, \mathcal{Z}} &= -\int_0^T \int_\Omega \big( \partial_t g_t + v\cdot\nabla_x g_t + S(f_t)\cdot\nabla_v g_t \big)f_t \dd z \dd t + \int_0^T (\frac{\dd}{\dd t}\dog - {\bf w})\cdot \varphi_d \dd t \\
 &\hspace*{2em}+ \int_\Omega g(T)f(T) - g(0)f(0) \dd z - \int_\Omega (f(0) -f_0) \eta_f \dd z   + (\dog(0)-\dog_0)\cdot \eta_d,
\end{align*}
with the adjoint state $(\varphi,\eta)\in \mathcal{X}\times \left( L^2(\Omega) \times \mathbb{R}^{MD} \right)$, $\varphi=(g,\varphi_d)$ and $\eta=(\eta_f,\eta_d)$. Similar to the microscopic case we define the extended Lagrangian corresponding to Problem \ref{Opt_PDE} as
\begin{equation*}
 \mathcal{L}_\infty(p,{\bf w},\varphi,\eta) = J_\infty(p,{\bf w}) + \langle e_\infty(p,{\bf w}),(\varphi,\eta) \rangle_{\mathcal{Z}^*, \mathcal{Z}},
\end{equation*}

Again, we denote the three parts of the cost functional by $J^i$, $i=1,2,3$, as
\begin{gather*}
 J_\infty^1(p) = \frac{\sigma_1}{4T} \int_0^T |\mathbb{V}_\infty(f_t) - \bar{V}_\infty|^2 \dd t, \qquad J_\infty^2(p) =  \frac{\sigma_2}{2T} \int_0^T |\mathbb{E}_\infty(f_t)-\Edes|^2 \dd t, \\
  J_\infty^3({\bf w})=\frac{\sigma_3}{2MT} \int_0^T \norm{\textbf{w}(t)}_{\mathbb{R}^{MD}}^2 \dd t.
\end{gather*}

Analogous to the microscopic case we derive the adjoint system and the optimality condition by calculating the derivatives of $\mathcal{L}_\infty$ w.r.t.~the state variable $p$ in direction $h^p=(h^f,h^d)\in\mathcal{Y}$, and the control $w$ in direction $h^w\in U$. The standard $L_2$-calculus yields
\begin{align*}
 d_f J_\infty^1 (p) [h^f] 
 &=\frac{\sigma_1}{T} \int_0^T \int_\Omega \Big(\mathbb{V}_\infty(f_t) - \bar{V}_\infty\Big) |x - \mathbb{E}_\infty(f_t)|^2  \dd h_t^f \dd t, \\
 d_f J_\infty^2 (p) [h^f] 
 &= \frac{\sigma_2}{T} \int_0^T \int_\Omega x\cdot (\mathbb{E}_\infty(f_t) - \Edes) \dd h_t^f \dd t, \\
 d_{{\bf w}} J_\infty^3 ({\bf w}) [h^w] &= \frac{\sigma_3}{MT} \int_0^T {\bf w}(t) \cdot h^w(t) \dd t.
\end{align*}
Let $\varphi = (g, \varphi_d)$ be the adjoint state corresponding to $p=(f,\dog)$. Then we obtain for the constraint part of the extended Lagrangian the following G\^ateaux derivatives:
\begin{align*}
 \langle d_f e_\infty(p,{\bf w})[h^f] , (\varphi,\eta) \rangle  &=  -\int_0^T \int_\Omega  \big( \partial_t g_t + v\cdot\nabla_x g_t + S(f_t)\cdot\nabla_v g_t \big)h_t^f \dd z \dd t, \\
 &\hspace*{4em}- \int_0^T \int_\Omega d_fS(f_t)[h_t^f]\cdot\nabla_vg_t \, f_t\dd z \dd t \\ 
 &\hspace*{4em}+ \int_\Omega g(T) h^f(T) -  h^f(0) g(0) - h^f(0) \eta_f \dd z,\\
 \langle d_{\dog} e_\infty(p,{\bf w})[h^d] , (\varphi,\eta)  \rangle &= \int_0^T \frac{\dd}{\dd t}h^d(t)\cdot \varphi_d(t) \dd t + h^d(0)\cdot \eta_d, \\
 &\hspace*{4em}- \int_0^T \int_\Omega d_{\dog} S(f_t)[h^d(t)]\cdot \nabla_vg_t \, f_t\dd z \dd t, \\
 \langle d_{{\bf w}}  e_\infty(p,{\bf w})[h^w] , (\varphi,\eta)  \rangle  &= - \int_0^T h^w(t) \cdot \varphi_d(t) \dd t.
\end{align*}
%
Assuming again the adjoint state $\varphi_d$ to be sufficiently regular, we may integrate by parts to obtain a strong formulation of the adjoint system. For the nontrivial terms we calculate  the following representations:
\begin{align}
 \int_\Omega d_fS(f_t)[h_t^f]\cdot\nabla_vg_t \, f_t\dd z
 &= -\int_\Omega \int_\Omega K_1(x,\bar x)\, h_t^f \dd \bar z \cdot \nabla_v g_t(z)\,f_t \dd z  \nonumber \\
 &= \int_\Omega \left(\int_\Omega K_1(\bar x,x)  \cdot \nabla_v g_t(z)\,f_t \dd z\right) h_t^f \dd \bar z \nonumber \\
 &=: \int_\Omega D_f(f_t)[g_t](z)\,h_t^f\dd z, \nonumber \\
 \int_\Omega d_{\dog} S(f_t)[h^d(t)]\cdot \nabla_vg_t \, f_t\dd z &= \int_\Omega \left( d_\dog S(f_t)[h^d(t)]\right) \cdot \nabla_v g_t(z)\,f_t \dd z \nonumber \\
 &= \left(\int_\Omega d_\dog S(f_t)[\nabla_vg_t(z)]\,f_t \dd z\right)\cdot h^d(t)\nonumber \\
 \label{D_fS}
 &=: D_\dog (f_t)[g_t]\cdot h^d(t). 
\end{align}
This yields the following adjoint system and optimality condition.

\begin{prop}\label{PDEKKT}
The optimality condition corresponding to Problem \ref{Opt_PDE} reads
 \begin{equation}\label{eq:var_vlasov}
  \int_0^T \left(\frac{\sigma_3}{MT} {\bf w}_*(t) - \varphi_d(t) \right)\cdot({\bf w}(t)-{\bf w}_*(t)) \dd t\ge 0  \qquad \text{for all\; ${\bf w} \in \mathcal{U}_{ad}$},
 \end{equation}
 where $\varphi=(g,\varphi_d)\in \mathcal{Y}$ satisfies the adjoint system given by
 \begin{subequations}\label{adSys}
 \begin{align}
  \partial_t g_t + v\cdot \nabla_x g_t &= -S(f_t) \cdot \nabla_v g_t + D_f(f_t)[g_t] - d_{\emph \x} J_\infty(t)  \label{ad_PDE}, \\
  \frac{\dd}{\dd t}\varphi_d &= -D_{\emph\dog}(f_t) [g_t], \label{ad_dog}
 \end{align}
 where
 \begin{equation}\label{dxJinf}
 d_{\emph \x} J_\infty(t) = \frac{\sigma_1}{T} \Big( \mathbb{V}(f_t) - \bar{V}_\infty \Big) |x-\mathbb{E}_\infty(f_t)|^2 + \frac{\sigma_2}{T} (\mathbb{E}_\infty(f_t) - \Edes) \cdot x \qquad\text{on\; $\Omega$}.
\end{equation}
 supplemented with the terminal conditions $g_T = 0$ and $\varphi_d(T) =  0$.

\end{subequations}
\end{prop}

\begin{rem}
 Note, that the optimality conditions of Problem \ref{Opt_ODE} and \ref{Opt_PDE} coincide. This is due to the fact that $J^3$ is identical in both problems. The mean-field limit only affects the adjoint system, but formally they can be identified in the following way: along the characteristics of the partial differential equation for $g$, $\xi_1$ corresponds to $\frac{1}{N}\nabla_x g$ and $\xi_2$ to $\frac{1}{N}\nabla_v g$. Finally, $\varphi_d$ can be identified directly with $\xi_3$.
\end{rem}

\begin{rem}
 The variational inequalities \eqref{eq:var_ode} and \eqref{eq:var_vlasov} may be equivalently expressed as fixed point problems in terms of a projection operator $\text{Proj}_U\colon U\to \mathcal{U}_{ad}$ which is defined by \cite{Pinnau}
 \[
  \text{Proj}_U(h)=\argmin\nolimits_{u \in \mathcal{U}_{ad}} \|u-h\|_U\qquad\text{for any\; $h\in U$}.
 \]
 Indeed, since $U$ is a Hilbert space the following statements are equivalent for any $k\in U$ and $\gamma>0$:
 \begin{enumerate}
  \item $u_*\in\mathcal{U}_{ad}$,\quad $\langle k,u-u_*\rangle_U\ge 0$\quad for all $u\in\mathcal{U}_{ad}$,
  \item $u_* = \text{Proj}_U(u_*-\gamma k)$.
 \end{enumerate}
 Consequently, the variational inequalities \eqref{eq:var_ode} and \eqref{eq:var_vlasov} may be expressed as
 \begin{align*}
  u_* = \text{Proj}_U(u_*-\gamma k)\in\mathcal{U}_{ad},
 \end{align*}
 where $k(\con)=\sigma_3/(MT) \con_* - \xi_3$ for the microscopic case \eqref{eq:var_ode} and $k({\bf w})=\sigma_3/(MT) {\bf w}_* - \varphi_d$ for the mean-field case \eqref{eq:var_vlasov}. In our particular case, $\text{Proj}_U$ has the explicit representation given by
 \begin{equation}\label{projection}
  \text{Proj}_U(h)(t) = \begin{cases}
  u_{\text{max}}\frac{h_m(t)}{|h_m(t)|} & \text{for\; $|h_m(t)|>u_{\text{max}}$}, \\
  h_m(t) & \text{otherwise},
  \end{cases}\quad m=1,\ldots,M,\quad \text{a.e. in\; $(0,T)$}.
 \end{equation}
\end{rem}

Due to the high dimensionality of the problem, it is not recommended to solve the system consisting of state equation, adjoint equation and optimality condition all at once. Therefore, we employ a gradient descent method. In order to define this method we make use of the gradient of the reduced cost functional which is derived in the following for both the microscopic problem and its mean-field counterpart.

\subsection{Gradient of the Reduced Cost Functional} \label{sec:redGrad}
In this section we introduce the reduced cost functionals $\hat{J}_N(\con)$ and $\hat{J}_\infty({\bf w})$ and formally calculate their gradients $\nabla \hat{J}_N(\con)$ and $\nabla \hat{J}_\infty({\bf w})$ which we need for the descent algorithms. 

By means of the control-to-state operators $\mathcal{G}_N\colon U\to Y$ and $\mathcal{G}_\infty\colon U\to\mathcal{Y}$ introduced in Section~\ref{sec:state}, we define the reduced cost functionals as
\begin{align*}
 \hat{J}_N(\con) := J_N(\mathcal{G}_N(\con),\con),\qquad
 \hat{J}_\infty({\bf w}) := J_\infty(\mathcal{G}_\infty({\bf w}),{\bf w}).
\end{align*}
Assuming sufficient regularity for $\mathcal{G}^N$ and $\mathcal{G}^\infty$ we further derive the gradients of the reduced cost functionals. Making use of the state equation $e_N(y,\con) = 0$  and $e_\infty(p,{\bf w})=0$ we implicitly obtain $d\mathcal{G}_N(\con)$ and $d\mathcal{G}_\infty({\bf w})$ via
\begin{align*}
0= d_\con e_N(\mathcal{G}_N(\con),\con) &= d_y e(\mathcal{G}_N(\con),\con)[d\mathcal{G}_N(\con)] + d_\con e_N(\mathcal{G}_N(\con),\con), \\ 
0= d_{\bf w} e_\infty(\mathcal{G}_\infty({\bf w}),{\bf w}) &= d_p e(\mathcal{G}_\infty({\bf w}),{\bf w})[d\mathcal{G}_\infty({\bf w})] + d_{\bf w} e_\infty(\mathcal{G}_\infty({\bf w}),{\bf w}),
\end{align*}
With the help of the adjoint equations
\begin{equation*}
 (d_ye(y,\con))^*[\xi] = - d_y J_N(y,\con)   \qquad \text{and} \qquad (d_pe(p,{\bf w}))^*[\varphi] = - d_p J_\infty(p,{\bf w})
\end{equation*}
we may calculate the G\^ateaux derivative of $\hat{J}$ in the direction $h\in U$, which gives
\begin{align*}
 d\hat{J}_N(\con)[h] &= \langle d_y J_N(y,\con), d\mathcal{G}_N(\con)[h]\rangle_{Y^*,Y} + \langle d_\con J_N(y,\con),h\rangle_U = \langle \frac{\sigma_3}{MT} \con -\xi_3, h \rangle_{U},\\
 d\hat{J}_\infty({\bf w})[h] &= \langle d_p J_\infty(p,{\bf w}), d\mathcal{G}_\infty({\bf w})[h]\rangle_{\mathcal{Y}^*,\mathcal{Y}} + \langle d_{\bf w} J_\infty(p,{\bf w}),h\rangle_U =\langle \frac{\sigma_3}{MT} {\bf w} -\varphi_d,h  \rangle_{U}.
\end{align*}
Since $U$ is a Hilbert space, we may use the Riesz representation theorem to identify the gradients
\begin{equation}\label{gradient}
 \nabla \hat{J}_N (\con) = \frac{\sigma_3}{MT} \con -\xi_3 \qquad \text{and} \qquad \nabla \hat{J}_\infty ({\bf w}) = \frac{\sigma_3}{MT} {\bf w} - \varphi_d.
\end{equation}
Now, we have all ingredients at hand to state the gradient descent method for the numerical simulations.

\section{Numerical Schemes and Algorithms}\label{sec:Alg}
Next, we investigate two control strategies, namely the Instantaneous Control (IC) and Optimal Control (OC). Both strategies are based on the adjoints and make use of the gradients derived in the previous section. As standard iterative adjoint based optimization methods, the algorithms rely on solvers for the state and adjoint systems and update the controls using the gradient information. 

\subsection{Numerics for the Forward and Adjoint IVP}
All occurring ODE systems are solved with the explicit fourth order Runge--Kutta solver \texttt{dopri} of the Python package \texttt{python.scipy}. The high order is more relevant for the adjoint problem than for the forward problem, since the adjoint system is stiffer due to the additional term arising from the linearization of the non-linear interaction and the cost functional. In order to have compatible data that can be passed from the forward problem to the adjoint problem and then further to the gradient update, we fix a time discretization in advance and give the solvers the freedom to make intermediate steps.

To obtain a numerical scheme that is independent of the number of particles involved we rescale the adjoint ODE by multiplying with $N$. This has the effect that the $N$-dependence of the terms in \eqref{terminal_ODE} emerging from the cost functional, is vanishing. In fact, for $i=1,\dots,N$, we set $r_i(t) = N \xi_i^1(t)$ and $s_i(t) = N \xi_i^2(t)$, and obtain the rescaled adjoint ODE system
\begin{subequations}\label{rescaled_ad}
\begin{align}
 \frac{\dd}{\dd t}r_i &= - \frac{1}{N} \sum_{j\ne i} \nabla _{x_i} K_1(x_i,x_j)(s_i-s_j) - \frac{1}{M} \sum_m\nabla_{x_i}K_2(x_i,d_m) s_i - \frac{1}{T}d_{x_i}J_N(t), \\
 \frac{\dd}{\dd t}s_i &= - r_i- \alpha s_i, \\
 \frac{\dd}{\dd t}\varphi_i &= \frac{1}{NM} \sum_{i=1}^N\nabla_{x_i} K_2(x_i,d_m) s_i. 
\end{align}
where
\begin{equation}
d_{x_i} J_N(t) = \sigma_1(\mathbb{V}(\x(t)) - \bar{V}_N)\big(x_i(t) - \mathbb{E}_N({\x}(t))\big) + \sigma_2\big(\mathbb{E}_N({ \x}(t))-\Edes\big),
\end{equation}
with terminal conditions $r(T) = 0$, $s(T)=0$ and $\varphi(T)=0$.
\end{subequations}

\subsection{Numerics for the Mean-field Equation and its Adjoint}
The forward and backward solves for the mean-field optimization are realized using the {\em Strang splitting} scheme\cite{StrangSplitting}. This scheme applies a semi-Lagrangian solver in spatial and a semi-implicit finite-volume scheme in velocity space. Using the following short hand notation for  \eqref{PDEsys}
\begin{equation*}
\partial_t f = -v \cdot \nabla_x f - \nabla_v  \cdot (S(f)  f),
\end{equation*}
we define the splitting
\begin{subequations}\label{discretization}
\begin{align}\label{dis_1}
\partial_t f^* &= -\frac{1}{2} \nabla_v \cdot ( S(f^*)   f^*) , &f^*(t) &= f(t), \\
\label{dis_2} \partial_t f^{**} &= - v \cdot \nabla_x f^{**}, &f^{**}(t) &=f^*(t+\tau), \\
\label{dis_3} \partial_t f &= -\frac{1}{2} \nabla_v  \cdot (S(f) f) , &f(t) &= f^{**}(t+\tau).
\end{align}
\end{subequations}
A Semi-Lagrangian method \cite{sonnendrucker1999semi,klar2009semi} is used to solve \eqref{dis_2}. The polynomial reconstruction  is done with the help of cubic Bezier curves. As usual, information is traced back along the characteristic curves starting in the respective grid point and interpolated at the end of the characteristic to get the values at the current time.
The discretizations in velocity space, \eqref{dis_1} and \eqref{dis_3}, use a second order finite volume scheme where the advection is approximated by a Lax-Wendroff flux \cite{leveque2002finite,quarteroni2008numerical}. Oscillations caused by non-smooth solutions are intercepted with the help of a van-Leer limiter \cite{VanLeerLimiter}. For more details on the scheme see \cite{Andreas}. 

Basically the same code is used for the adjoint system, for this sake rewrite
\begin{equation*}
S(f) \cdot \nabla_v g = -\nabla_v \cdot (S(f) g) + 2 \alpha g.
\end{equation*}
Further, we need to add the term resulting from the linearization of the non-linear interaction and the cost functional. Altogether we obtain for the adjoint system the splitting
\begin{align*}
\partial_t g^* &= -d_x J_\infty + \frac{1}{2}  \left( -\nabla_v \cdot (S(f) g^*) + 2 \alpha g^* + D_f(f_t)[g_t]\right) , &g^*(t) &= g(t), \\
 \partial_t g^{**} &= - v \cdot \nabla_x g^{**}, &g^{**}(t) &=g^*(t+\tau), \\
\partial_t g &= \frac{1}{2} \left( -\nabla_v \cdot (S(f) g) + 2 \alpha g  + D_f(f_t)[g_t] \right), &g(t) &= g^{**}(t+\tau).
\end{align*}
with $D_f(f_t)[g_t]$ as defined in \eqref{D_fS} and $d_x J_\infty$ as in \eqref{dxJinf}.

\subsection{The Optimization Algorithms}
In the following we discuss the algorithms for the IC and OC approach employed for the optimization. In general the procedures differ in the amount of information used to compute the gradient to update the control for the next iteration. The descent directions are based on the gradients in the following manner. Let $k$ denote the current iteration, $s_k$ the descent direction and $c_k$ the control $\con_k$ or $\textbf{w}_k$, respectively. Then the control is updated due to
\begin{equation}\label{eq:update}
 \tilde{c}_{k+1} = c_k + \omega_k s_k,
\end{equation}
where $\omega_k$ is an admissible step size.
\begin{rem}
Note that the interpretation of one iteration differs for the IC ansatz and the OC approach. For IC the iteration $k$ denotes the current time slice whereas for OC $k$ is the counter of the optimization iterations.
\end{rem}
For the Instantaneous Control method we use steepest descent steps to update the control once on every time slice, therefore $s_k = -\nabla \hat{J}(u_k)$. For the Optimal Control approach we apply a non-linear conjugate gradient method (NCG), see Algorithm \ref{CGdirection} for details. In the pseudocode the gradient $\nabla \hat{J}(u_k)$ is denoted by $q_k$ to shorten the notation. We choose NCG, since it converges faster than the steepest descent method and requires less memory than BFGS methods. Nevertheless, the CG factor $\gamma_{k}$ is chosen such that the descent direction resulting from Algorithm \ref{CGdirection} equals the descent direction of a BFGS formula applied one time to the identity \cite[Thm 5.13]{numopt}.
\RestyleAlgo{boxruled}
\begin{algorithm2e}
\caption{Non-Linear Conjugate Gradient Direction}\label{CGdirection}
 \KwData{current gradient $q_k$, previous gradient $q_{k-1}$, previous descent direction $s_{k-1}$}
 \KwResult{descent direction $s_k$}
 initialization\;
 \eIf{ k = 1 }{ $s_k = -q_k$}{
  $\gamma_{k-1} = \frac{(q_k-q_{k-1},q_k)}{(q_k-q_{k-1},s_{k-1})}$ \qquad \text{and} \qquad  $s_k = -q_k - \gamma_{k-1} s_{k-1}$\;   
\If{ $(s_k,q_k) > -\tol_{CG}$}{ $s_k = -q_k$}}
\end{algorithm2e}
Due to the maximal velocity constraint \eqref{eq:control} we need to make sure that computed the controls are feasible. We therefore project the controls onto the feasible set in every iteration using the operator $\text{Proj}_U$ defined in \eqref{projection}. The step sizes $\omega_k$ are obtained with the help of the projected Armijo step size rule \cite{Pinnau} (see Algorithm \ref{ArmRule}).
\RestyleAlgo{boxruled}
\begin{algorithm2e}[h!]
\caption{Projected Armijo Stepsize Rule }\label{ArmRule}
 \KwData{Current control $c_k$, gradient $q_k$, initial $\omega_0$, initial $\gamma$}
 \KwResult{new control $c_{k+1}$}
 initialization\;
 \While{ $\hat{J}( \emph{Proj}_U[c_k - \omega_k q_k]) \ge  \hat{J} + \gamma\omega_k \norm{q_k}^2$ }{
   $\omega_k = \omega_k/2$
 }
\end{algorithm2e}

The Instantaneous Control is a feedback control strategy. The gradient is based on state information of the current time slice. Figuratively, that means that the external agents may react instantaneously on the behavior of the crowd and do not take the past and future into account. Thus, the adjoints and gradients are based on the forward data of one time step only, resulting in advantageous memory consumption (cf. Table~\ref{memory} below) and short computation times. 
To realize the IC approach we divide the time interval $[0,T]$ into $K$ equidistant slices and sequentially compute the optimal control for every time slice (cf. \cite{HinzeVolkwein}). The control problem  then reads \\[3mm]
\textit{ For $k = 1,\dots,K$ find $u_k \in \mathbb{R}^{DM}$ such that 
\begin{align*}
 u_k = \argmin_{u_k \in \mathcal{U}_{ad}}  \hat{J}_N(u_k) &= \argmin_{u_k\in \mathcal{U}_{ad}} \frac{\sigma_1}{4} (\mathbb{V}_N({\emph \x}_k) - \mathbb{V}_N({\emph \x}_0))^2 + \frac{\sigma_2}{2} \norm{\mathbb{E}_N({\emph \x}_k)-\Edes}_2^2 + \frac{\sigma_3}{2M} \norm{u_{k}}_{2}^2  \\
&\text{subject to}\qquad  y_k = y_{k-1} + dt \cdot F(y_{k-1},u_k)
\end{align*}}
in the ODE case and
\\[3mm]
\textit{ For $k = 1,\dots,K$ find $w_k \in \mathbb{R}^{DM}$ such that 
\begin{align*}
 w_k = \argmin_{w_k \in \mathcal{U}_{ad}}  \hat{J}_\infty(w_k) &= \argmin_{w_k \in \mathcal{U}_{ad}} \frac{\sigma_1}{4} (\mathbb{V}_\infty(f_k) - \mathbb{V}_\infty(f_0))^2 + \frac{\sigma_2}{2} \norm{\mathbb{E}_\infty(f_k)-\Edes}_2^2 + \frac{\sigma_3}{2M}  \norm{w_{k}}_{2}^2  \\
&\text{subject to} \qquad p_k = p_{k-1} + dt \cdot G(p_{k-1},w_k)
\end{align*}}
in the PDE case. 

The control for time slice $k+1$ is initialized with
\begin{equation}\label{IC_update_control}
	c_{k+1} = 0.1 c_k
\end{equation}
to make sure that the penalty term of the cost functional is initially non-zero.
The final positions and velocities of particles and the positions of the external agents of the time step $t_{k-1}$ are the initial values for the time slice $t_k$. The controls are assumed to be constant in each time slice. The controls for the whole simulation interval $[0,T]$ are obtained by gluing the solutions on the time slices together leading to a piecewise constant control function. A pseudocode for the IC procedure can be found in Algorithm \ref{ICAlg}. 
\RestyleAlgo{boxruled}
\begin{algorithm2e}[h!]
\caption{Instantaneous Control Algorithm}\label{ICAlg}
 \KwData{Initial data of $y$ or $p$; parameter values}
 \KwResult{Instantaneous control $\con$ or $\textbf{w}$; the corresponding states $y(\con)$ or $p(\textbf{w})$; optimal functional values}
 initialization\;
 $t=0; dt=T/K$\;
 \While{$t < T$}{
 \begin{itemize}
 \item [0.] solve state system \eqref{ODEsys} or \eqref{PDEsys}\;
 \item [1.] solve adjoint problem given in \eqref{rescaled_ad} or \eqref{adSys}\;
 \item [2.] compute gradient corresponding to \eqref{gradient}\;
 \item [3.] compute step size using the Armijo rule \eqref{ArmRule}
 \item [4.] update controls by steepest descent step \eqref{eq:update}\;
 \item [5.] project control into the feasible set using \eqref{projection}\;
 \item [6.] initialize controls for the next time slice corresponding to \eqref{IC_update_control}\;
 \end{itemize}
  t = t+dt
 }
\end{algorithm2e}
\begin{rem}
 Note that we perform only one gradient step on each time slice. Numerical tests showed that more gradient steps marginally improve the result but drastically increase the computational time.
\end{rem}

The Optimal Control strategy is an open-loop approach that works on the whole time interval $[0,T]$ at once. The adjoint and gradient computations use the information of a full forward solve on the interval $[0,T]$. This results in a large amount of memory storage needed (cf. Table~\ref{memory} below).  A vivid description of this strategy is that the agents plan their whole journey in advance. 
The optimization algorithm operates as follows. First, the state system is solved for initially given velocities for the external agents. Its solution is then used to compute the adjoint solution backward in time on the full interval $[0,T]$ as well. Finally, the control for Problem \ref{Opt_ODE} or Problem \ref{Opt_PDE} is updated for the next iteration using the gradient corresponding to \eqref{gradient} and \eqref{eq:update}. The iterations are stopped when the $L_2$ norm of two consecutive controls is smaller than the tolerance $\tol$ see \eqref{stopping}. A pseudocode for the Optimal Control Algorithms is given in Algorithm \ref{OCAlg}. 
\RestyleAlgo{boxruled}
\begin{algorithm2e}
\caption{Optimal Control Algorithm}\label{OCAlg}
 \KwData{Initial data of $y$ or $p$; parameter values}
 \KwResult{Optimal control $\con$ or $\textbf{w}$; the corresponding states $y(\con)$ or $p(\textbf{w})$; optimal functional values}
 initialization\;
 solve state system \eqref{ODEsys} or \eqref{PDEsys}\;
 \While{$eps_{rel} >  \tol$}{
 \begin{itemize}
 \item [1.] solve adjoint problem given in \eqref{rescaled_ad} or \eqref{adSys}\;
 \item [2.] compute gradient corresponding to \eqref{gradient}\;
 \item [3.] update controls by NCG step \eqref{eq:update} with descent direction of Algorithm \ref{CGdirection}\;
 \item [4.] project control onto the feasible set using \eqref{projection}\;
 \item [5.] find appropriate step size by projected Armijo-rule using Algorithm \ref{ArmRule}\;
 \item [6.] compute $\norm{\nabla \hat{J}}$, update $eps_{rel}$ \;
 \end{itemize}
 }
\end{algorithm2e}

\begin{rem}
Note, that the line search using projected Armijo Rule might require multiple solves of the forward system. 
\end{rem}

\section{Numerical Results}\label{sec:NumRes}
The numerical simulations for the mean-field equation are conducted on the domain $\Omega \subset \mathbb{R}^{4}$ as $ \Omega = [-100,100]^2 \times [-5,5]^2$, i.e., we set $D=2$. For the mean-field simulation we use the scaled CFL condition
\begin{equation}
 \frac{\tau |V|T}{Lh} \le 0.5
\end{equation}
with $L = 200, |V| = 5$. The grid parameter $h$ is varied throughout the simulations to investigate the convergence of the scheme. In fact, we use $25,50$ or $100$ grid points in each of the two directions leading to $h = 0.04, 0.02$ or $0.01$.  Our particular choice of the interaction potentials are the Morse potentials as proposed in  \cite{Dorsogna, CarrilloMatrinPanferov}. For fixed positive parameters  $A_j,a_j,R_j,r_j$ we have 
\begin{equation}\label{InteractionPotentials}
 \Phi_j(x-y) = R_j \exp\left(-\frac{|x-y|}{r_j}\right) - A_j \exp\left( - \frac{|x-y|}{a_j}\right), \quad j=1,2.
\end{equation}
The parameters $A_j,R_j$ denote the attraction and repulsion strengths and $a_j, r_j$ the radius of interaction. The case $j=1$ refers to the interaction of the individuals, the interaction of individuals with external agents is denoted by $j=2$. Inspired by \cite{NoScaling}, we use the values
\begin{equation*}
 A_1 = 20,\quad R_1=50,\quad a_1=100,\quad r_1=2, \qquad\qquad A_2 = 5,\quad R_2=100,\quad a_2=1000,\quad r_2=50.
\end{equation*}

\begin{rem}
Note that the coefficients of the interaction potentials are scaled according to $\Omega$. As mentioned in the introduction we had the model problem of dogs herding sheep to the destination $Z$ in mind. The parameters for the sheep-sheep interaction have long range attraction and repulsion on a very short range. The parameters modeling the sheep-dog interaction have a larger repulsive influence in order to reflect the guiding property correctly. 
\end{rem}

\begin{figure}[htbp]
	\begin{minipage}{0.49\textwidth} 
	\includegraphics[width=1.\textwidth]{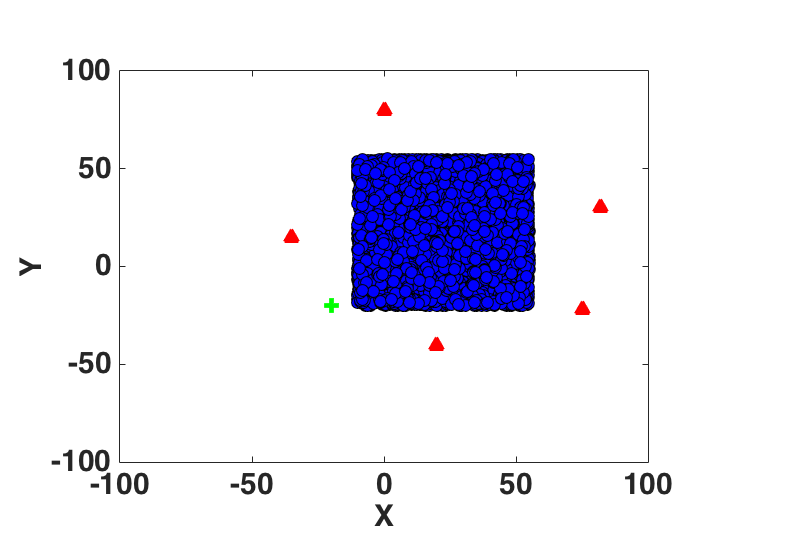}
	\end{minipage}
	\hfill
	\begin{minipage}{0.49\textwidth}
	\includegraphics[width=1.\textwidth]{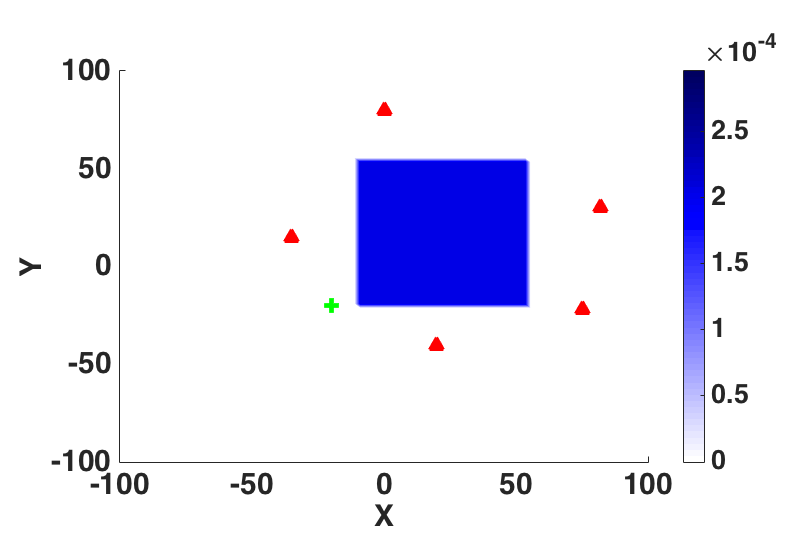}
	\end{minipage}
\caption{Initial configurations. Left: microscopic case with 8000 sheep represented by blue dots. Right: initial probability distribution $f$ for the mean-field case. In both cases the positions of the dogs $d$ are marked by red triangles and the destination $\Edes$ by a green dot.}
\label{fig:initial}
\end{figure}
Further parameters are fixed as follows. The destination $\Edes = (-20,-20)$ is depicted in green. The initial support of the sheep is set to $\Omega_0 = [-10,55]\times[-20,55]$, i.e., $f_0(x,v)$ is the uniform distribution with support $\Omega$ and the initial positions and velocities of the sheep are chosen by realizations of i.i.d. random variables with $\text{law}(f_0)$. In the Figures the sheep and dogs are represented by blue markers and red triangles, respectively. The initial configurations of the microscopic and the mean-field case are shown in Figure~\ref{fig:initial}. The parameters for the following results are set to
\begin{equation*}
 T=10,\quad \bar{V}_n = 0.9\mathbb{V}_N(\x_0),\quad \bar{V}_\infty = 0.9\mathbb{V}_\infty(f_0) \quad \text{and} \quad \sigma_3 = 0.000001.
\end{equation*}
Thus, the desired variance is $10\%$ less than the initial variance given by the initial distribution of the crowd. Since the values of the cost functional descent are in the range of $10e-6$ on each time slice of the IC Algorithm, we have to choose $\sigma_3$ rather small. Otherwise, the Armijo line search does not succeed.
The Armijo parameters for IC and OC are set to
\begin{equation*}
 \omega_0^{{\emph IC}} = 1000, \qquad \omega_0^{{\emph OC}} = 10.
\end{equation*}
\begin{rem}
Note, that in the given time interval $[0,T]$ the task of steering the crowd to the destination $\Edes$ is impossible to realize. In fact, our main focus lies on the comparison of the behaviour of the (sub-)optimal controls as $N$ increases. Here sub-optimal refers to the IC simulations and optimal to the OC results.
\end{rem}

\subsection{Numerical Results Using the IC Algorithm}
\begin{figure}
 	\begin{minipage}{0.49\textwidth}
	\includegraphics[width=1.\textwidth]{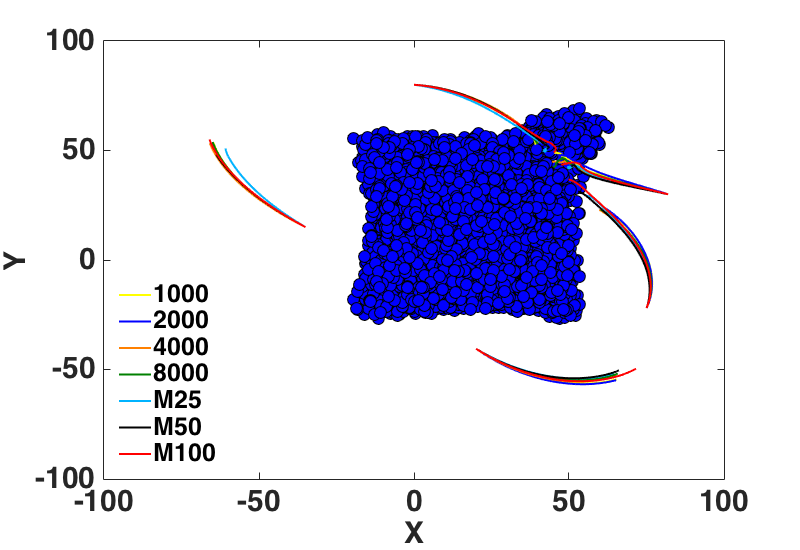}
	\end{minipage}
	\hfill
	\begin{minipage}{0.49\textwidth}
	\includegraphics[width=1.\textwidth]{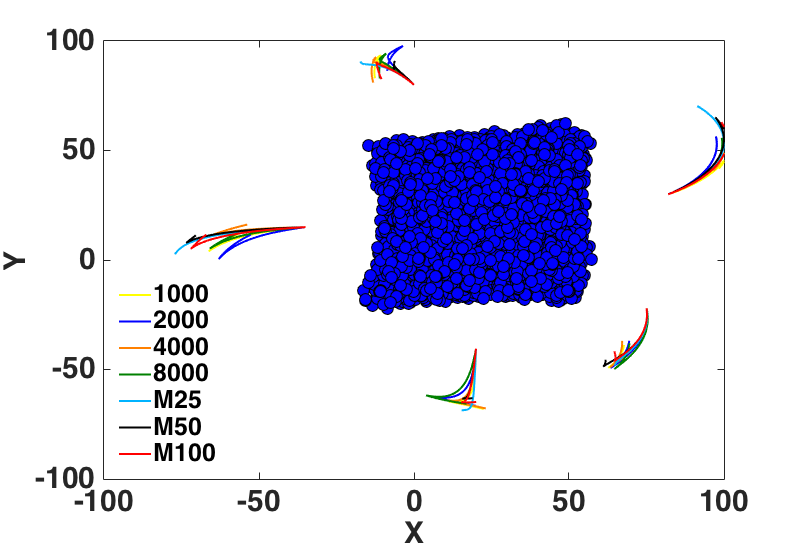}
	\end{minipage}
	\caption{Influence of the cost functional parameters. Left: Setting IC\ref{S1} - the dogs positioned at the left side and the bottom of the crowd clear the way to the destination. The others push the crowd from behind towards the destination. Right: Setting IC\ref{S2} - the dogs keep distance and orient towards to corners to reduce the variance of the crowd. These observations perfectly agree with the intention we had when modeling the cost functional.}
	 \label{fig:influence}
\end{figure}

In this section the results obtained by the IC algorithm are illustrated. First, we discuss the influence of the cost functional parameters. We therefore set up three test cases: \ref{S1} stresses the variance part of the cost functional, we expect the dogs to move towards the corners of the crowd in order to reduce the spread of the crowd. In $\ref{S2}$ the second part of $J$ is emphasized, thus the dogs shall push the crowd towards the destination. For test case \ref{S3} we choose the weights of the cost functional such that the focus lies on steering the crowd to the destination $\Edes$ while the variance term has minor influence but cannot be neglected. In particular, we use the following parameter values
\begin{align}\label{S1} \tag{\textbf{S1}}
 \sigma_1 &= 0.09, &&\sigma_2 = 0.001, \\
 \label{S2} \tag{\textbf{S2}}
 \sigma_1 &= 0.0001, &&\sigma_2 = 0.9, \\
 \label{S3} \tag{\textbf{S3}}
 \sigma_1 &= 0.005, &&\sigma_2 = 0.5. 
\end{align}
These choices assure that $J_1$ of test case \ref{S1} has the same order of magnitude as $J_2$ of test case \ref{S2} and the other way around. The following notation is used in the Figures
\begin{equation*}
 J_1 = \frac{\sigma_1}{4} \Big(\mathbb{V}_N(\x_t) - \bar{V}_N\Big)^2 \quad\text{and}\quad J_2 = \frac{\sigma_2}{2} \norm{\mathbb{E}_N(\x_t)-\Edes}_2^2
\end{equation*}
for the microscopic cost functional and analogous for the mean-field cost functional. The spatial and velocity grid is discretized with the same amount of grid points $25, 50$ or $100$ in each of the two directions for the IC algorithm. Due to the large amount of memory needed for the OC algorithm, we use grids of $25$ and $50$ points in this case. Further information about the memory consumption of the different simulations can be found in Table~\ref{memory}. The graphs corresponding to mean-field solutions are labelled M$\#$ where $\#$ denotes the number of grid points. The graphs corresponding to microscopic simulations are denoted by their respective number of particles. In Figure~\ref{fig:influence} the influence of the cost functional parameters is illustrated. 

\begin{figure}
	\begin{minipage}{0.49\textwidth}
	\includegraphics[width=1.\textwidth]{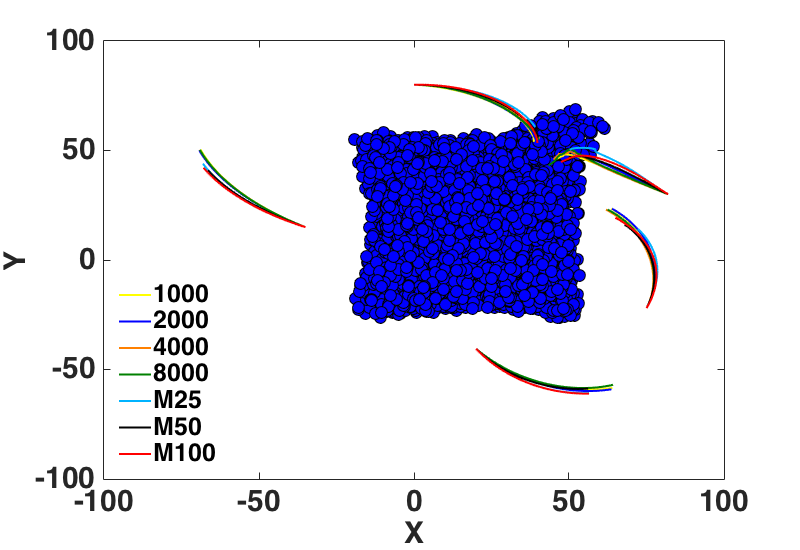}
	\end{minipage}
	\hfill
	\begin{minipage}{0.49\textwidth}
	\includegraphics[width=1.\textwidth]{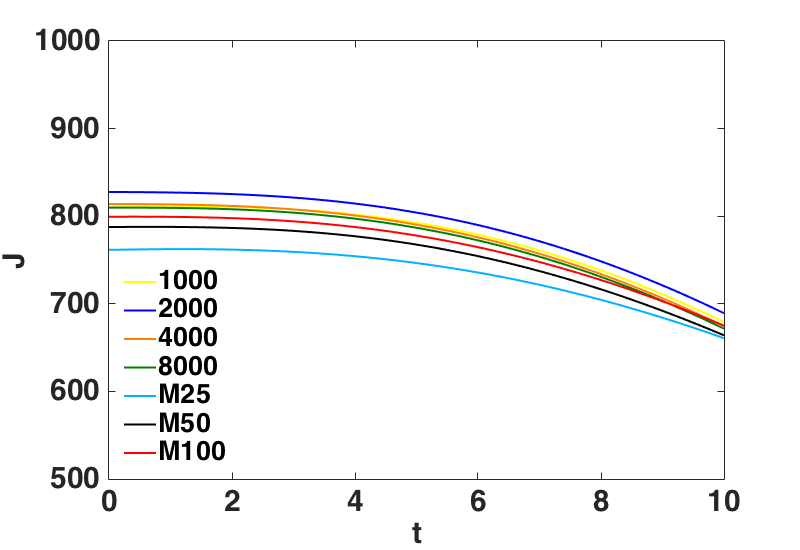}
	\end{minipage}
	\caption{Setting IC\ref{S3} - Trajectories of the dogs and evolution of the cost functional. Left:  the dogs positioned at the left side and at the bottom of the crowd clear the way to the destination. The others push from behind. Note that the dogs on the left to not enter the crowd as deep as in setting IC\ref{S1}. Right: Evolution of the cost functional values.}
	\label{fig:pathAndCost}
\end{figure}

\begin{figure}
	\begin{minipage}{0.49\textwidth}
	\includegraphics[width=1.\textwidth]{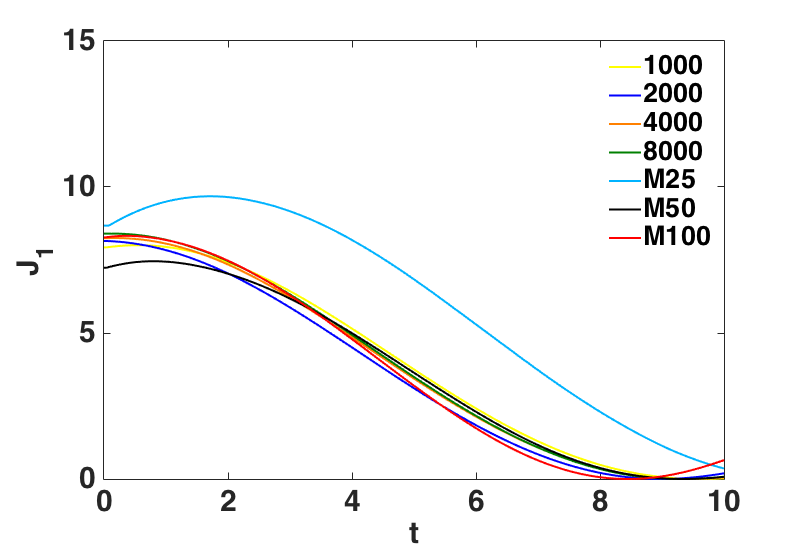}
	\end{minipage}
	\hfill
	\begin{minipage}{0.49\textwidth}
	\includegraphics[width=1.\textwidth]{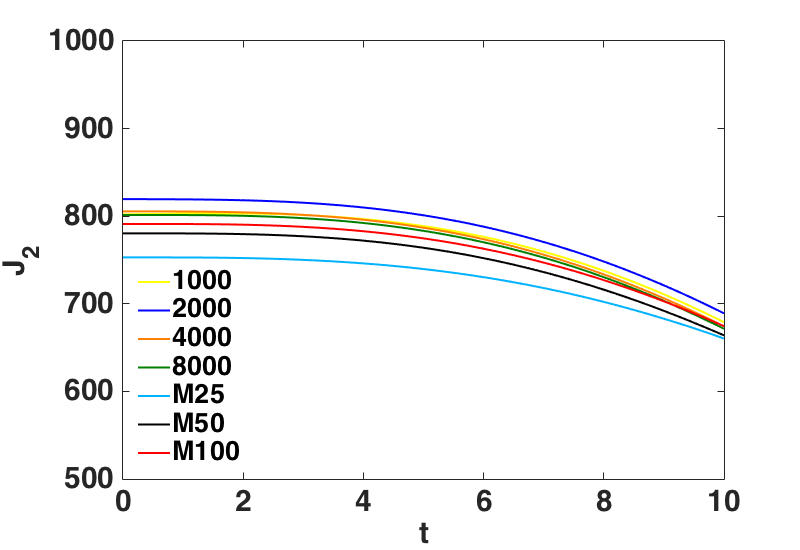}
	\end{minipage}
	 \caption{Setting IC\ref{S3} - Evolution of the cost functional parts. Left: values of $J_1$. The simulation $M25$ overestimates the variance significantly. All other graphs are in good agreement. Right: values of $J_2$. All graphs show similar behavior. We see convergence as $N\rightarrow \infty$ and as $h \rightarrow 0$.}
	 \label{fig:costDetail}
\end{figure}

We assume that \ref{S3} is the most realistic setting for our purpose of leading the crowd to the destination, since the focus lies on $J_2$ which measures the distance to the destination while the influence of $J_1$ is not negligible. In Figure~\ref{fig:pathAndCost}(left) the trajectories of the dogs and the crowd at $T=10$ are shown. We note that the dogs are not entering the crowd as deep as in test case \ref{S1}. 
Comparing the values of $J$, $J_1$ and $J_2$ in Figure~\ref{fig:pathAndCost}(right) and Figure~\ref{fig:costDetail} we obviously find that the part measuring the crowd's distance to the destination is superimposing the others. Note that the graphs of the cost functional values are given with respect to time, i.e. it may happen that the values increase. Since the initial positions of the dogs is chosen arbitrarily, they move to appropriate positions first, causing a slight increase of the cost functional.  Afterwards the cost functional is decreasing as expected. The simulation on the coarse mean-field discretization $M25$ overestimates the variance term significantly as can be seen in Figure~\ref{fig:costDetail}(left). The evolution of $J_1$ is in good resemblance for all other discretization. Similarly, there is accordance in the graphs showing the evolution of $J_2$. Furthermore, we see the convergence of the graphs to $M100$ as $N$ increases and the convergence of the mean-field simulations as the grid is refined.

In Table~\ref{Table_IC} we investigate the convergence for $N \rightarrow \infty$ using the following scaled norms
\begin{subequations} \label{num_norms}
\begin{align}
\norm{u - u_\text{ref}} &= \frac{1}{TMV} \int_0^T \norm{u(t) - u_\text{ref}(t)}_{ \mathbb{R}^4} \dd t, \\
\norm{J- J_\text{ref}}&= \frac{1}{T} \int_0^T |J(t) -J_\text{ref}(t) |  \dd t, \\
\norm{\rho^N - \rho_\text{ref}} &=\frac{1}{T} \int \norm{\rho(t,\tilde{x},\tilde{y}) - \rho_\text{ref}(t,\tilde{x},\tilde{y})} \dd \tilde{x} \dd \tilde{y} \dd t,
\end{align}
\end{subequations}
where $\dd \tilde{x}$ and $\dd \tilde{y}$ are the spacial measures in x and y direction rescaled by $1/L$. As reference values we use the results of the simulation $M100$. To compute the norms the empirical density $\rho^N$ of the microscopic simulations is approximated by a histogram which is based on the grid of the corresponding mean-field simulation. 
In the first two columns we see the convergence of the mean-field scheme, as expected the values are decreasing as the grid is refined. The four columns on the right illustrate the behavior as $N$ increases. All norm values are decreasing for increasing number of particles from $N=1000$ to $N=8000$. The convergence is more transparent for the integral quantities as for the velocities (cf. Remark~\ref{integral_quantities}).

\begin{rem}\label{integral_quantities}
Since mean-field quantities are averaged quantities it is very common to compare integral values like $J$ and $\rho^N$ in this setting. Note that the velocities and the trajectories of the agents are no such quantities. We therefore expect to have blurred values for $\norm{ u - u_\text{ref}}$. 
\end{rem}

\begin{table}[h!]
\begin{tabular}{| c | c | c || c | c | c | c | }
 \hline
  IC\ref{S1}       &  M25    & M50    & 1000     & 2000  & 4000  & 8000         \\
 \hline
 $\norm{J - J_\text{ref}}$ &    4.8820e-02 & 7.6180e-03  &  3.4355e-02   & 4.2821e-02  & 3.2420e-02 & 3.0061e-02   \\
 \hline
 $\norm{ u- u_\text{ref}}$ &  7.8537e-02   & 2.8145e-02  & 2.4854e-02  & 2.4726e-02   &   2.4609e-02  & 2.4733e-02     \\
  \hline
 $ \norm{\rho^N - \rho_{\text{ref}}}$ &  - & - & 1.2422e-05   & 1.2396e-05   & 1.2288e-05  &    1.22281e-05       \\
  \hline
  IC\ref{S2}     &  M25    & M50    & 1000     & 2000  & 4000  & 8000         \\
 \hline
 $\norm{J - J_\text{ref}}$ &  3.3725e-02 &  5.231e-03 & 5.442e-03  &  1.756e-03   & 1.908e-03  & 1.919e-03     \\
 \hline
 $\norm{ u- u_\text{ref}}$ &  6.3913e-02   & 1.9453e-02  & 2.3938e-02  & 2.5209e-02   &   2.4165e-02 & 2.4628e-02     \\
  \hline
 $ \norm{\rho^N - \rho_{\text{ref}}}$ &  - & - & 1.2590e-05   & 1.2581e-05   & 1.2460e-05  &    1.2454e-05      \\
  \hline
   IC\ref{S3}      &  M25    & M50    & 1000     & 2000  & 4000  & 8000         \\
 \hline
 $\norm{J - J_\text{ref}}$ &  2.7053e-02 &  4.826e-03 & 5.510e-03  &  1.1599e-02   & 5.050e-03  & 3.579e-03     \\
 \hline
 $\norm{ u- u_\text{ref}}$ &  5.2202e-02   & 1.8787e-02  & 2.2655e-02  & 2.2637e-02   &   2.2717e-02 & 2.26106e-02     \\
  \hline
 $ \norm{\rho^N - \rho_{\text{ref}}}$ &  - & - & 1.2426e-05   & 1.2404e-05   & 1.2292e-05  &    1.2298e-05      \\
  \hline
\end{tabular}
 \caption{Investigation of convergence as $N\rightarrow \infty$ using the norms in \eqref{num_norms} with the respective values of $M100$ as reference values.}
 \label{Table_IC}
\end{table}

\begin{figure}
 	\begin{minipage}{0.49\textwidth}
	\includegraphics[width=1.\textwidth]{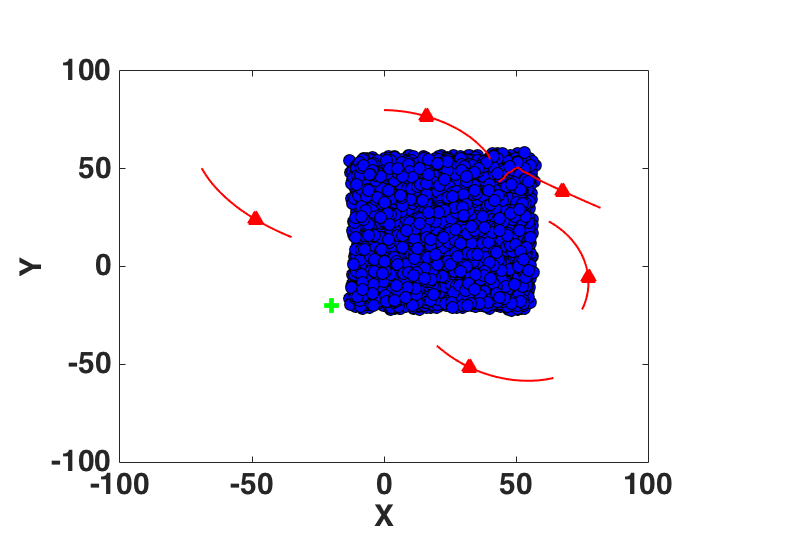}
	\end{minipage}
	\hfill
	\begin{minipage}{0.49\textwidth}
	\includegraphics[width=1.\textwidth]{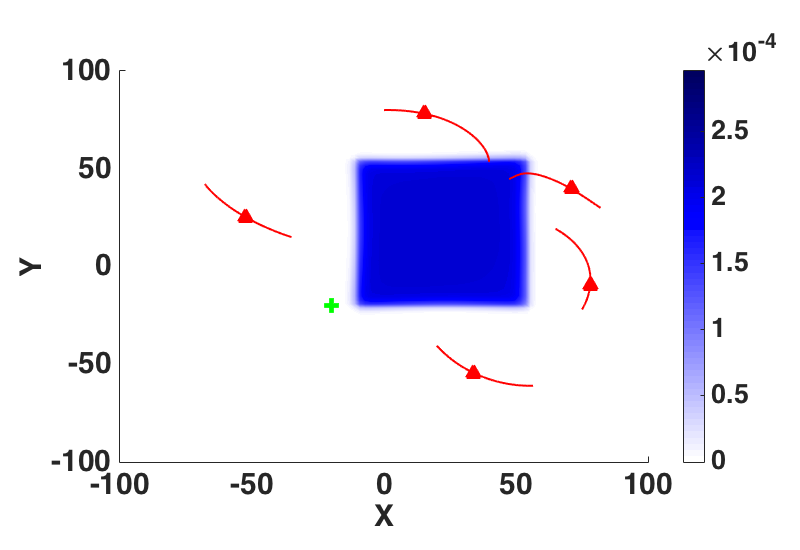}
	\end{minipage}
	\begin{center} Microscopic and mean-field crowd together with the dogs at $t=3.3$. \end{center}
 	\begin{minipage}{0.49\textwidth}
	\includegraphics[width=1.\textwidth]{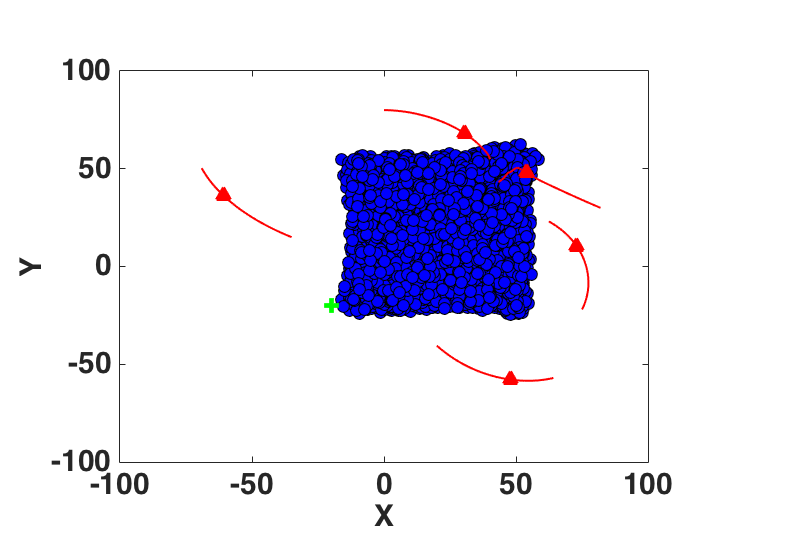}
	\end{minipage}
	\hfill
	\begin{minipage}{0.49\textwidth}
	\includegraphics[width=1.\textwidth]{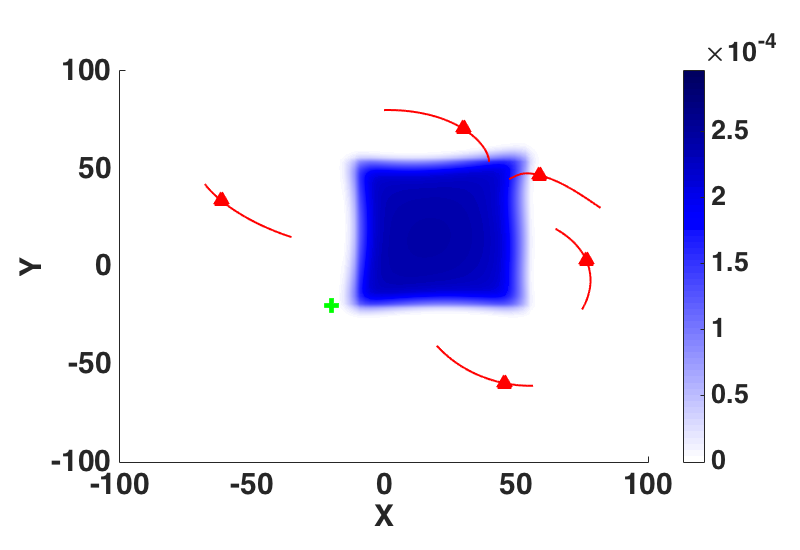}
	\end{minipage}
	\begin{center} Microscopic and mean-field crowd together with the dogs at $t=6.6$. \end{center}
 	\begin{minipage}{0.49\textwidth}
	\includegraphics[width=1.\textwidth]{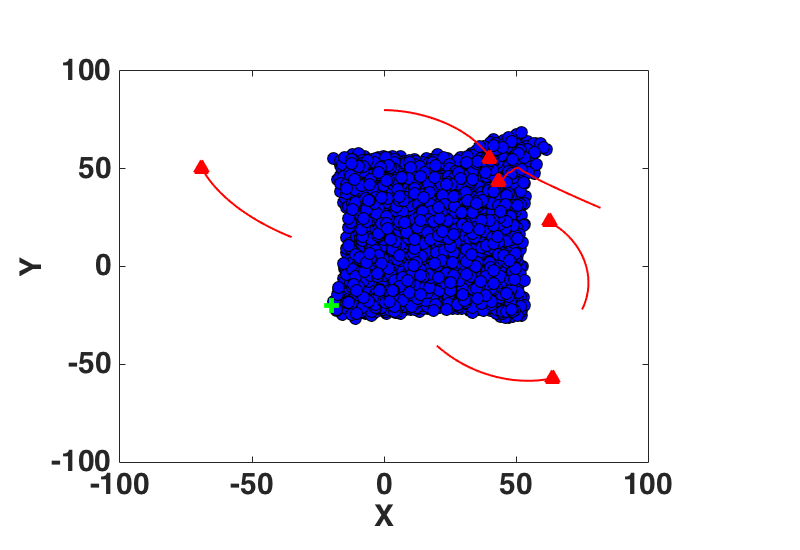}
	\end{minipage}
	\hfill
	\begin{minipage}{0.49\textwidth}
	\includegraphics[width=1.\textwidth]{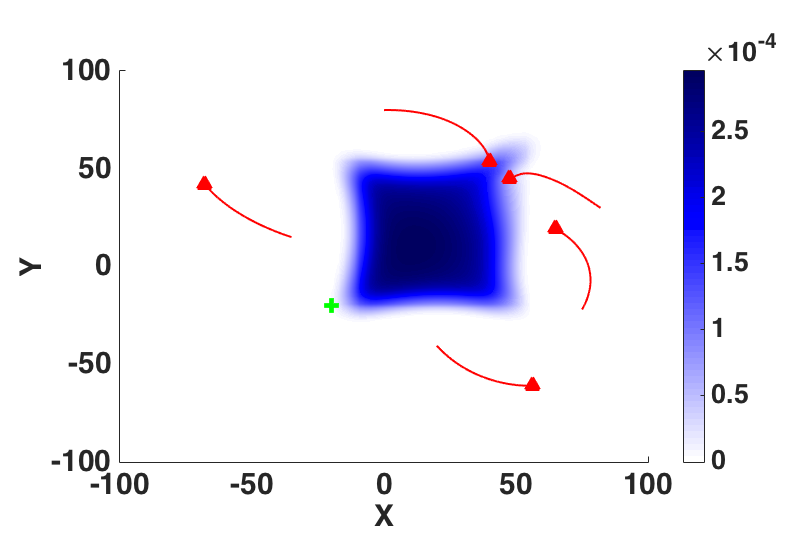}
	\end{minipage}
	\begin{center} Microscopic and mean-field crowd together with the dogs at $t=10$. \end{center}
	\caption{Setting IC\ref{S3} - Comparison of the microscopic $(N=8000)$ and mean-field (M100) crowd and the dogs at different times computed with the IC simulation.}
	\label{fig:ParticleMFVergleich}
\end{figure}

In Figure~\ref{fig:ParticleMFVergleich} the evolution of the crowd and the dogs is illustrated for the microscopic ($N=8000$) and the mean-field case $M100$. The red lines show the trajectories of the dogs as above. Additionally to the information in Figure~\ref{fig:pathAndCost}, the red markers indicate the current positions of the dogs at different times. These positions and the form of the crowd are in good accordance.

\subsection{Numerical Results Using the Optimal Control Approach}
The results for the OC approach are discussed in this section. As this approach is of open-loop type,  we stop the iteration, if 
\begin{equation}\label{stopping}
 \frac{\norm{w_{k+1} - w_0}_{L^2((0,T),\mathbb{R}^{MD})}}{ \norm{w_{0}}_{L^2((0,T),\mathbb{R}^{MD})}} \le \tol,
\end{equation}
with $\tol = 0.05$ for the computations. Due to the huge amount of memory needed for the simulations we study the case \textbf{S3} for $N=1000,2000,4000$ and  $M = 25, 50$ only. %
\begin{figure}[h!]
 	\begin{minipage}{0.49\textwidth}
	\includegraphics[width=1.\textwidth]{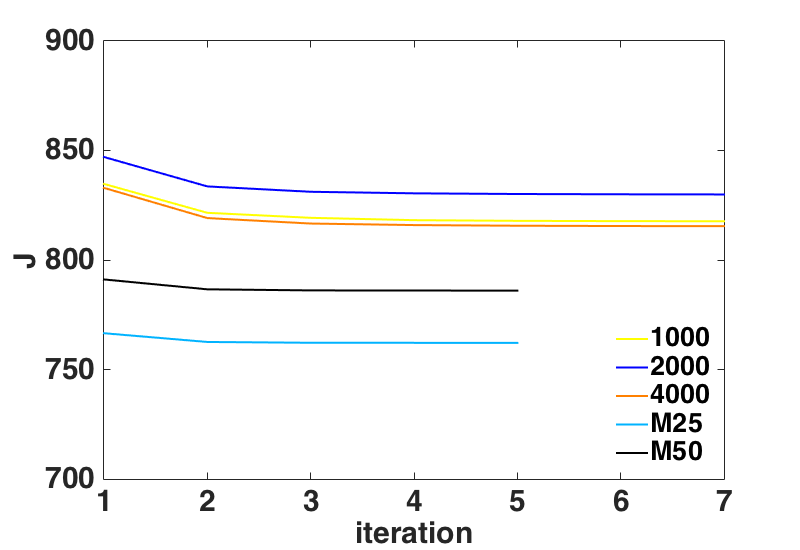}
	\end{minipage}
	\hfill
	\begin{minipage}{0.49\textwidth}
	\includegraphics[width=1.\textwidth]{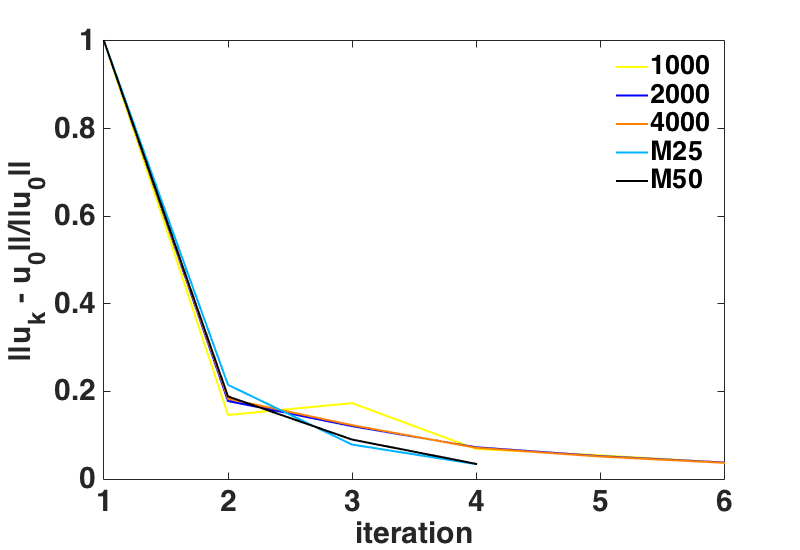}
	\end{minipage}
	 \caption{Setting OC\ref{S3} - Left: The values of $J$ plotted over the OC iterations for different $N$ in the particle case and different grid sizes in the mean-field case. Right: The values of the relative norms used to evaluate the stopping criterion plotted over the iterations for different $N$ in the particle case and different grid sizes in the mean-field case.}
	 \label{fig:JandEpsOC}
\end{figure}

\begin{rem}
Note, that the reference values of the cost functional graphs in Figure~\ref{fig:pathAndCost}(right), Figure~\ref{fig:costDetail} and Figure~\ref{fig:JandEpsOC} differ. For the IC simulations the graphs are plotted over time. For the OC simulations the graph is plotted with respect to the iterations, the values have to decrease in this case. Of course, the values of $J_N$ and $J_\infty$ are in the same range.
\end{rem}

\begin{figure}
 	\begin{minipage}{0.49\textwidth}
	\includegraphics[width=1.\textwidth]{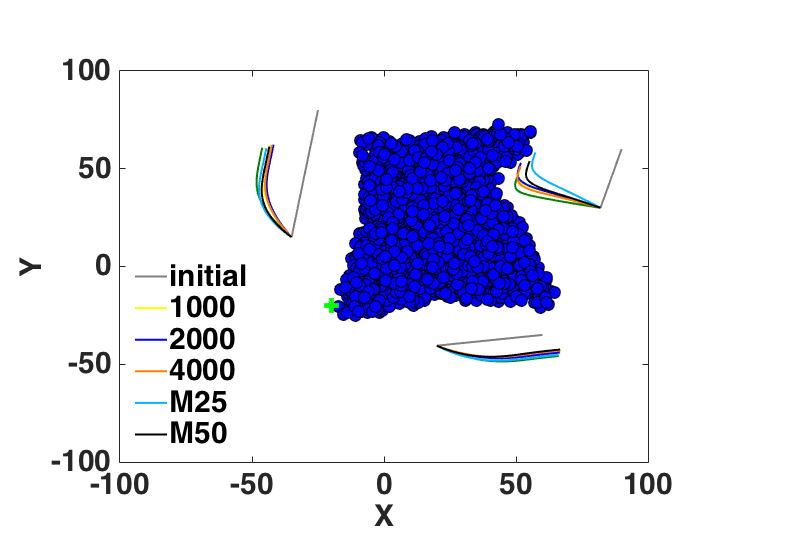}
	\end{minipage}
	\hfill
	\begin{minipage}{0.49\textwidth}
	\includegraphics[width=1.\textwidth]{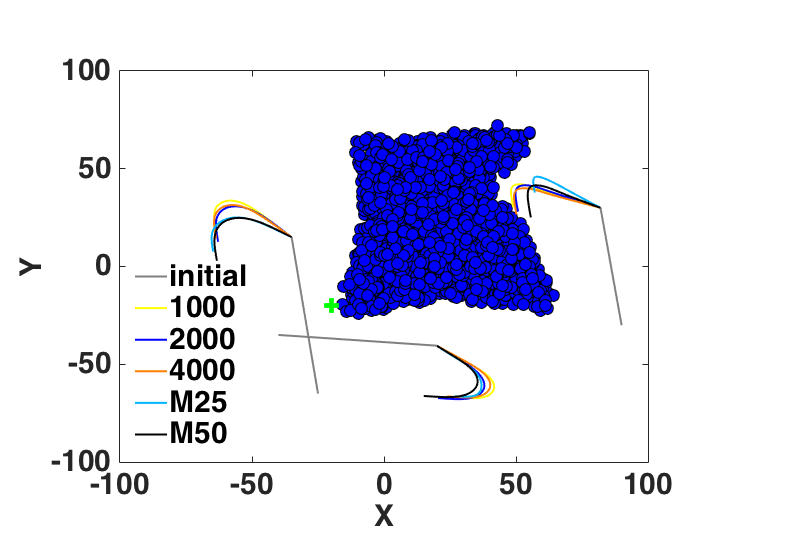}
	\end{minipage}
	 \caption{Setting OC\ref{S3} - Left:  Trajectories of the dogs for different $N$ in the microscopic case and different grid sizes in the mean-field case. Right:  Trajectories of the dogs for different initial controls. A comparison of the figures on the left and the right side shows the dependence of the solution on the initial data.}
	 \label{fig:PATHSandVelNormOC}
\end{figure}

\begin{figure}
 	\begin{minipage}{0.49\textwidth}
	\includegraphics[width=1.\textwidth]{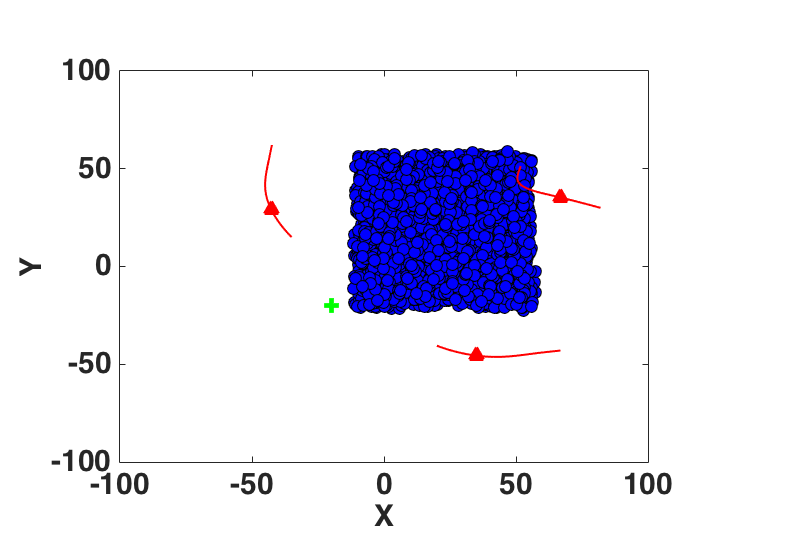}
	\end{minipage}
	\hfill
	\begin{minipage}{0.49\textwidth}
	\includegraphics[width=1\textwidth]{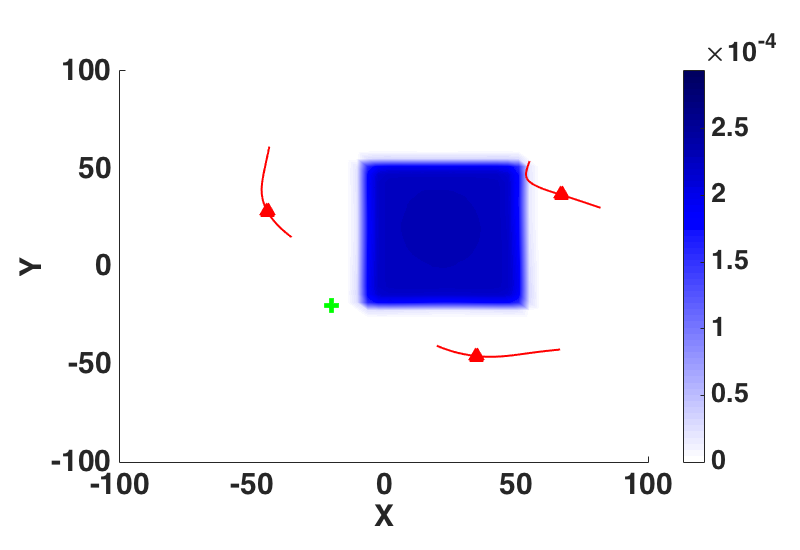}
	\end{minipage}
	\begin{center} Microscopic and mean-field crowd together with the dogs at $t=3.3$. \end{center}
 	\begin{minipage}{0.49\textwidth}
	\includegraphics[width=1.\textwidth]{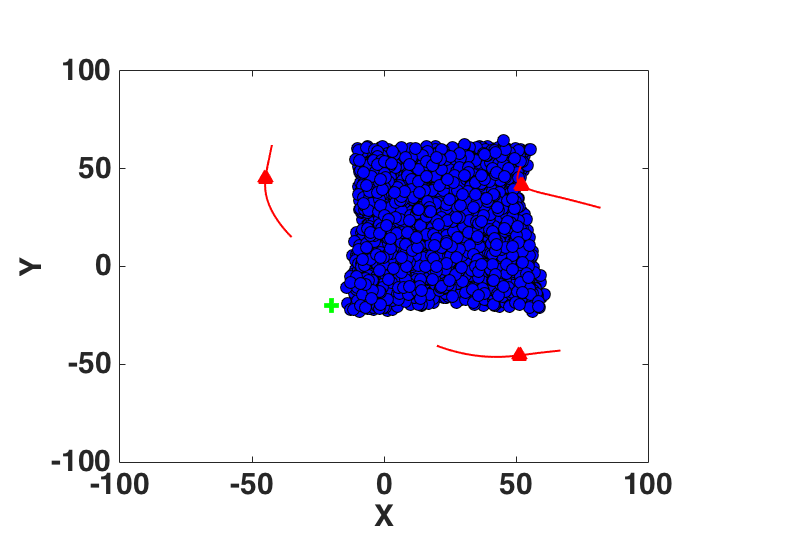}
	\end{minipage}
	\hfill
	\begin{minipage}{0.49\textwidth}
	\includegraphics[width=1.\textwidth]{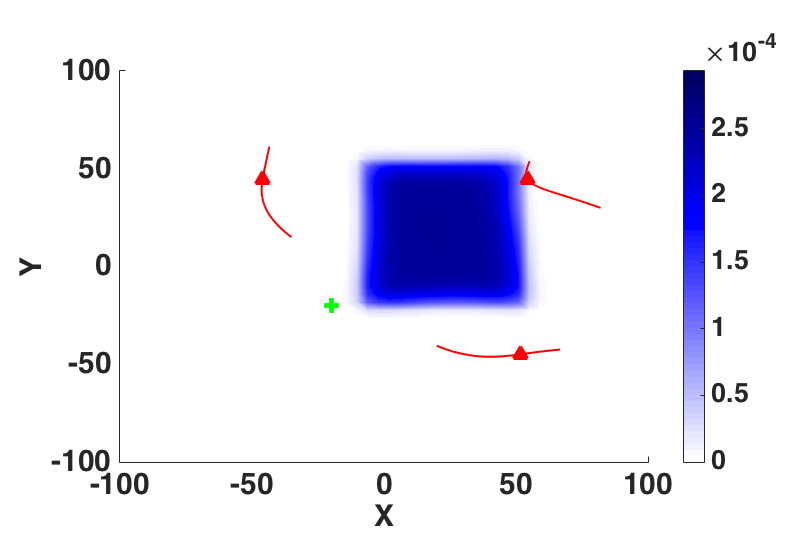}
	\end{minipage}
	\begin{center} Microscopic and mean-field crowd together with the dogs at $t=6.6$. \end{center}
 	\begin{minipage}{0.49\textwidth}
	\includegraphics[width=1.\textwidth]{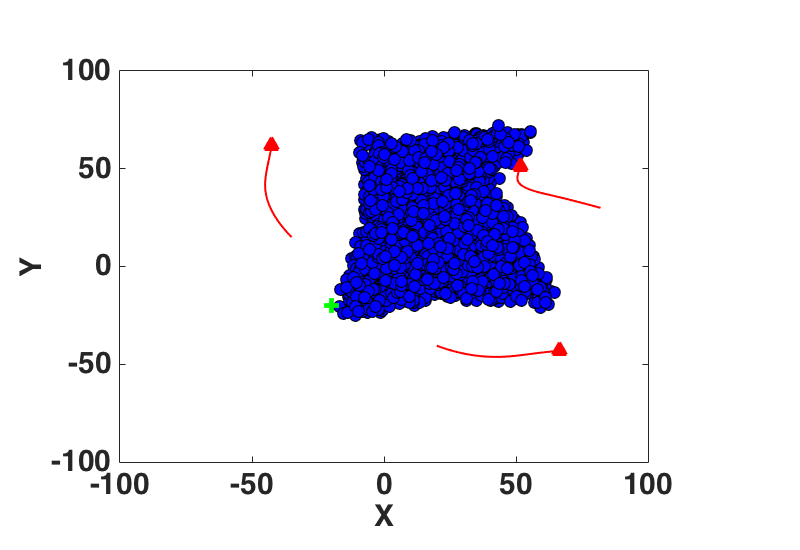}
	\end{minipage}
	\hfill
	\begin{minipage}{0.49\textwidth}
	\includegraphics[width=1.\textwidth]{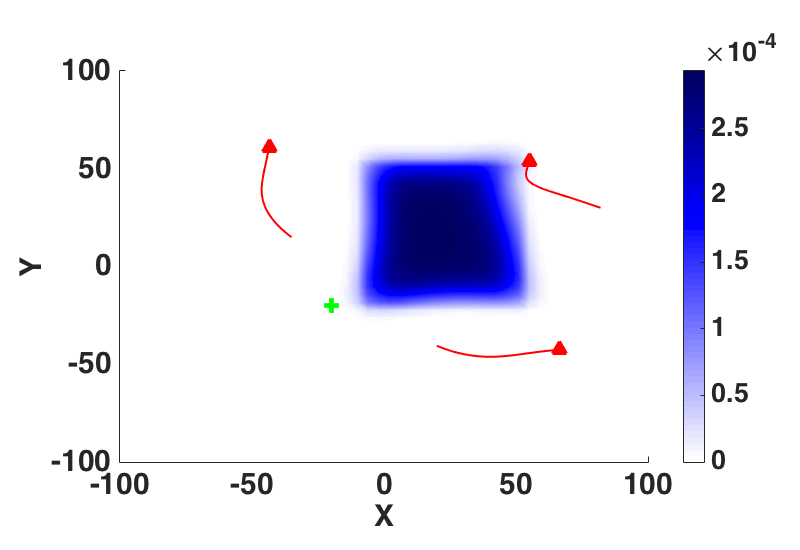}
	\end{minipage}
	\begin{center} Microscopic  and mean-field crowd together with the dogs at $t=10$. \end{center}
	\caption{Setting OC\ref{S3} - Comparison of the microscopic $(N=4000)$ and mean-field (M50) crowd and the dogs at different times computed with the OC Algorithm.}
	\label{fig:ParticleMFVergleichOC}
\end{figure}

The plots of the evolution of the cost functionals over the iterations in Figure~\ref{fig:JandEpsOC}(left) show the typical behavior of Optimal Control computations. In the first optimization iteration the most reduction is obtained followed by small improvements as the iterations proceed. As we have already discussed we cannot expect the dogs to steer the crowd to the destination in the short time interval given. Likewise, the cost functionals do reduce only marginal on this short interval. The figure on the right shows the relative norms which are used to evaluate the stopping criterion. Analogous to the behavior of the cost functionals, the graphs of the relative changes in velocity behave very similar for all microscopic and mean-field simulations. 
In Figure~\ref{fig:PATHSandVelNormOC} (left) we see the trajectories resulting from the optimal controls for the different number of particles and the mean-field simulations. The trajectories corresponding to the initial controls are depicted in gray. Obviously, in all test cases the deviation from the initial path due to the optimization is similar. Even though one gets the impression that the dog on the right of the mean-field simulation stays far away from the crowd in Figure~\ref{fig:PATHSandVelNormOC} (left), we observe in Figure~\ref{fig:ParticleMFVergleichOC} that this distance is caused by the resolution of the grid.  
\begin{table}
\begin{tabular}{| c | c | c |  }
	\hline
	  & $\norm{u-u_{\text{ref}}} $ & $\norm{\rho-\rho_{\text{ref}}} $  \\
	\hline
	1000 &   2.5311e-02 &  9.0509e-05  \\
	\hline
	2000 &   2.5353e-02  & 9.0347e-05  \\
	\hline
	4000 & 2.5208e-02 &   8.9968e-05 \\
	\hline
	M25   & 2.4585e-02   & - \\
	\hline
\end{tabular}
\caption{Setting OC\ref{S3} - Differences of the controls measured in the scaled norms given in \eqref{num_norms}.}
\label{tab:NormOC}
\end{table}
We investigate the convergence of the controls and the density of the crowd as $N \rightarrow \infty$ with the help of the norms defined in \eqref{num_norms} in Table~\ref{tab:NormOC}. For the microscopic simulations the density was computed using a histogram based on the grid of the mean-field simulation $M50$.

Further, the OC algorithm allows for an investigation of the dependence on the initial data. To study this dependence we describe different initial controls for the optimization procedure, which lead to the results shown in Figure~\ref{fig:PATHSandVelNormOC} (right). It is obvious that different initial data lead to different optimal controls. Nevertheless, the formation of the crowd is in good accordance.

In Figure~\ref{fig:ParticleMFVergleichOC} we see the evolution of the sheep crowd and the dogs at different times computed with the help of the OC algorithm. The red lines are the trajectories of the dogs. The red triangles denote their current position at the respective time $t$. The dogs positions and the form of the crowd is very similar at all time instances shown.

We end this section with the investigation of the memory consumed for the simulations. In Table~\ref{memory} we find that the the OC simulations need much more memory as the IC computations. Thus, future investigations of the optimization problem on a larger time interval can only be realized with the help of the IC approach.

\begin{table}
\begin{tabular}{| c | c | c | c | c | c | c | c | c  |}
 \hline
 N & 1000     & 2000  & 4000  & 8000  & grid & 25    & 50     & 100     \\
 \hline
 IC &  0.15   & 0.55  & 2.15  &  9.05 &  IC    &  0.08 & 0.977  &  15.473 \\
 OC &  0.17   & 0.6   & 3.21  &    -  &    OC  & 0.82  & 24.75  &  -      \\
  \hline
\end{tabular}
 \caption{Memory Consumption [Gigabyte]}
 \label{memory}
 \end{table}
 
\section{Conclusions}
We numerically verified for both approaches, IC and OC, that the controls of the microscopic problem converge to the corresponding mean-field controls as $N$ increases. From a computational point of view the IC approach has the advantage of being less memory consuming. In addition, in the case of dogs guiding a crowd of sheep, the IC ansatz is more realistic due to the feasibility of adjusting the velocity of the dogs in every time step and thus instantaneous reaction on the collective behaviour of the crowd. Nevertheless, one may think of applications, for example guides leading visitors through museums or exhibitions,  where it makes sense to plan journeys in advance and not react instantaneously, since the crowds may have clashed already. In this example the OC approach would be more convenient. Summarizing, our  simulations  indicate that the microscopic simulations for large number of particles can be appropriately approximated with the help of the mean-field simulations also in the control settings.

\section*{Acknowledgements} 
MB acknowledges support by ERC via Grant EU FP 7 - ERC Consolidator Grant 615216 LifeInverse.

\bibliographystyle{plain}
\bibliography{biblio}
\end{document}